\documentclass[10pt]{amsart}
\usepackage{amsmath,amscd,amssymb,epsfig}

\newcommand{\ZZ}{\mathbb{Z}}
\newcommand{\RR}{\mathbb{R}}
\newcommand{\SSS}{\mathbb{S}}

\newcommand{\sym}{\mathfrak{S}}
\newcommand{\alt}{\mathfrak{A}}

\newcommand{\inv}{\operatorname{inv}}
\newcommand{\ltrmin}{\operatorname{lrmin\,}}

\newcommand{\Des}{{\rm{Des}}}
\newcommand{\des}{{\rm{des}}}
\newcommand{\WDes}{{\rm{Nasc}}}
\newcommand{\HatWDes}{\widehat{\rm{Nasc}}}

\newcommand{\symmdiff}{\triangle}

\newcommand{\Pal}{{\tt{P}}}
\newcommand{\T}{{\tt{T}}}
\newtheorem{theorem}{Theorem}[subsection]
\newtheorem{corollary}[theorem]{Corollary}
\newtheorem{proposition}[theorem]{Proposition}

\newtheorem{defn}[theorem]{Definition}

\newtheorem{example}[theorem]{Example}
\newtheorem{remark}[theorem]{Remark}
\newtheorem{question}[theorem]{Question}

\numberwithin{figure}{section}

\begin{document}
\title{Alternating subgroups of Coxeter groups}

\author{Francesco Brenti}
\address{Dipartimento di Matematica, Universit\'{a} di Roma ``Tor
Vergata'', Via della Ricerca Scientifica, 00133 Roma, Italy}
\email{brenti@mat.uniroma2.it}

\author{Victor Reiner}
\address{School of Mathematics\\
University of Minnesota\\
Minneapolis, MN 55455\\
USA}
\email{reiner@math.umn.edu}

\author{Yuval Roichman}
\address{Department of Mathematics\\
Bar-Ilan University\\
52900 Ramat-Gan\\
Israel}
\email{yuvalr@math.biu.ac.il}



\thanks{
Second author supported by NSF grant DMS-0245379. Third author
supported in part by the Israel Science Foundation, founded by the
Israel Academy of Sciences and Humanities, grant nu. 947/04.}

\begin{abstract}
  We study combinatorial properties of the alternating subgroup of a Coxeter group,
using a presentation of it due to Bourbaki.
\end{abstract}

\date{February 5, '07}
\maketitle

\section{Introduction}
\label{intro-section}

For any Coxeter system $(W,S)$, its {\it alternating subgroup}
$W^+$ is the kernel of the {\it sign character} that sends every
$s \in S$ to $-1$. An exercise from Bourbaki gives a simple
presentation for $W^+$, after one chooses a generator $s_0 \in S$.  The
goal here is to explore
 the combinatorial properties of this presentation,
distinguishing in the four main sections
of the paper different levels of generality (defined below)
regarding the chosen generator $s_0$:
$$
\begin{matrix}
 & &s_0\text{ arbitrary}& & \\
 & &(Section~\ref{general-section})& & \\
 &\diagup&                    &\diagdown& \\
s_0\text{ evenly-laced }& & & &s_0\text{ a leaf }\\
(Section~\ref{evenly-laced-section})& & & &(Section~\ref{leaf-section})\\
 &\diagdown&                    &\diagup& \\
 & &s_0\text{ an even leaf }& & \\
 & &(Section~\ref{even-leaf-section})& &
\end{matrix}
$$

Section~\ref{general-section} reviews the presentation and explores some of
its consequences in general for the {\it length function},
{\it parabolic subgroups}, a {\it Coxeter-like complex} for $W^+$, and
the notion of {\it palindromes}, which play the role
usually played by {\it reflections} in a Coxeter system.  This section also defines
{\it weak} and {\it strong} partial orders on $W^+$ and poses some basic questions about them.

Section~\ref{evenly-laced-section} explores the
special case where $s_0$ is {\it evenly-laced}, meaning that
the order $m_{0i}$ of $s_0 s_i$ is even (or infinity) for all $i$.
It turns out that, surprisingly,
this case is much better-behaved. Here the unique, length-additive
factorization $W=W^J \cdot W_J$ for parabolic subgroups of $W$
induces similar unique length-additive factorizations within
$W^+$. One can compute generating functions for $W^+$ by
length, or jointly by length and certain descent statistics.  Here
the palindromes which shorten an element determine that element
uniquely, and satisfy a crucial {\it strong exchange property}.
This gives better characterizations of the weak and strong partial
orders, and answers affirmatively all the questions about these
orders from Section~\ref{general-section} in this case.

Section~\ref{leaf-section} examines how the general presentation simplifies to
what we call a {\it nearly Coxeter} presentation when $s_0$ is a {\it leaf} in
the Coxeter diagram, meaning that $s_0$ commutes with all but one of the other
generators in $S-\{s_0\}$.  Such leaf generators occur in many situations,
e.g. when $W$ is finite\footnote{Combinatorial aspects
of this nearly Coxeter presentation were explored for $W$ of type $A$ in
\cite{RegevRoichman}, and partly motivated the current work.}
and for most affine Weyl groups.

\begin{figure}
\epsfxsize=120mm
\epsfbox{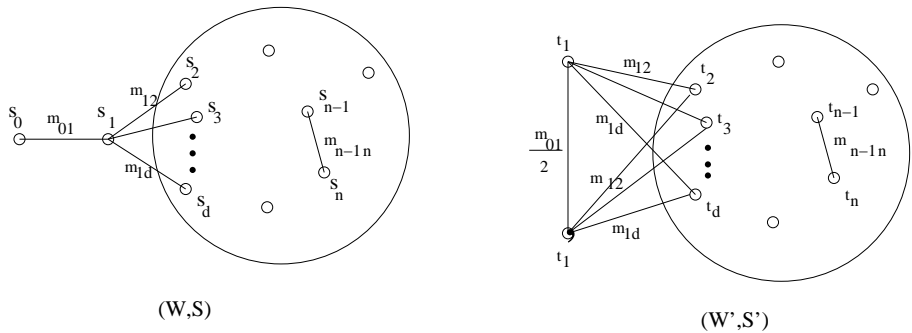}
\caption{Schematic of the relation between the diagrams for a Coxeter
system $(W,S)$ with even leaf node $s_0$, and the Coxeter system $(W',S')$ derived
from it, closely connected to the alternating group $W^+$.  The unique
neighbor of $s_0$ has been labelled $s_1$, so that $m_{01}$ is {\it even}.}
\label{Oriflamme-figure}
\end{figure}

Section~\ref{even-leaf-section} studies the further special case where $s_0$ is
an evenly-laced leaf.  The classification
of finite and affine Coxeter systems shows that {\it all} evenly-laced nodes $s_0$ are
even leaves when $W$ is finite, and this is almost always the case for $W$ affine.
In particular, even leaves occur in the finite type $B_n=(C_n)$ and the affine types
$\tilde{B_n}, \tilde{C_n}$.  When $s_0$ is an even leaf, there is an
amazingly close connection
between the alternating group $W^+$ and a {\it different} index $2$ subgroup $W'$,
namely the kernel of the homomorphism $\chi_0$ sending $s_0$ to $-1$ and all other Coxeter
generators to $+1$.
It turns out that this subgroup $W'$ is a (non-parabolic) reflection subgroup of $W$,
carrying its own Coxeter presentation $(W',S')$, closely related to the
Coxeter presentation of
$(W,S)$. This generalizes the inclusion of type $D_n$ inside $B_n$, and although
$W^+ \not\cong W'$, the connection allows one to
reduce all the various combinatorial questions for the presentation $(W^+,R)$
(length function, descent sets, partial orderings, reduced words) to their well-studied
counterparts in the Coxeter system $(W',S')$.


\tableofcontents

\section{The general case}
\label{general-section}

\subsection{Bourbaki's presentation}
\label{Bourbaki-subsection}

Let $(W,S)$ be a Coxeter system with generators
$S=\{s_0,s_1,\ldots,s_n\}$, that is, $W$ has a presentation of the
form
\begin{equation}
\label{W-presentation}
\begin{aligned}
W=\langle S=\{s_0,s_1,\ldots,s_n\} :
(s_i s_j)^{m_{ij}}=e \text{ for }0 \leq i \leq j \leq n \rangle
\end{aligned}
\end{equation}
where $m_{ij} =m_{ji} \in \{2,3,\ldots\} \cup \{\infty\}$ and $m_{ii}=2$.

The {\it sign character} $\epsilon: W \rightarrow \{\pm 1\}$ is
the homomorphism uniquely defined by $\epsilon(s) = -1$ for all $s
\in S$. Its kernel $W^+:=\ker(\epsilon)$ is an index two
subgroup called the {\it alternating subgroup} of $W$.

Once one has distinguished $s_0$ in $S$ by its zero subscript, an exercise in
Bourbaki \cite[ Chap. IV, Sec. 1, Exer. 9]{Bou} suggests a
simple presentation for $W^+$, which we recall here and prove along the
lines suggested by Bourbaki.

\begin{proposition}
\label{Bourbaki-exercise}
Given a Coxeter system $(W,S)$ with distinguished generator $s_0$,
map the set $R=\{r_1,\ldots,r_n\}_{i=1,2,\ldots,n}$ into $W^+$
via $r_i \mapsto s_0 s_i$.  Then this gives a set of generators
for $W^+$ with the following presentation:
\begin{equation}
\label{general-W-plus-presentation}
\begin{aligned}
&W^+ \cong \langle R=\{r_1,\ldots,r_n\} : \\
&\qquad \qquad r_i^{m_{0i}}=(r_i^{-1} r_j)^{m_{ij}}=e \text{ for }1\leq i<j \leq n \rangle.
\end{aligned}
\end{equation}
\end{proposition}
\begin{proof}
Consider the abstract group $H^+$ with the presentation by
generators $R$ given on the right side of
\eqref{general-W-plus-presentation}. One checks that the set map
$\alpha: R \rightarrow H^+$ sending $r_i$ to $r_i^{-1}$ extends to
an involutive group automorphism $\alpha$ on $H^+$:
the relation $(r_i r_j^{-1})^{m_{ij}}=e$ follows from the relation
$(r_i^{-1}  r_j)^{m_{ij}}=e$ in $H^+$ by taking the inverse of
both sides and then conjugating by $r_j$.

Thus the group $\ZZ/2\ZZ=\{1,\alpha\}$ acts on $H^+$, and one can form
the semidirect product $H^+ \rtimes \ZZ/2\ZZ$ in which
$(h_1\alpha^i) \cdot (h_2 \alpha^j )=h_1 \alpha^i(h_2) \cdot \alpha^{i+j}$.
This has either of the following two presentations:
$$
\begin{aligned}
&H^+ \rtimes \ZZ/2\ZZ \\
&\cong \langle r_1,\ldots,r_n, \alpha: \quad
                         \alpha^2=r_i^{m_{0i}}=(r_i^{-1} r_j)^{m_{ij}}=e \text{ for }1\leq i<j \leq n ,
                          \,\, \alpha r_i \alpha = r_i^{-1}  \rangle\\
&\cong \langle r_0, r_1,\ldots,r_n, \alpha:\\
       & \qquad \qquad r_0=\alpha^2=(r_i^{-1} r_j)^{m_{ij}}=e \text{ for }0 \leq i<j \leq n , \\
       & \qquad \qquad \alpha r_i \alpha = r_i^{-1}  \rangle\\
\end{aligned}
$$

We claim that the following two maps are well-defined and inverse isomorphisms:
$$
\begin{matrix}
W       & \overset{\rho}{\longrightarrow}  & H^+ \rtimes \ZZ/2\ZZ &\\
   s_i  & \longmapsto     & \alpha r_i (= r_i^{-1} \alpha) &\text{ for }i=1,\ldots,n\\
   s_0  & \longmapsto     & \alpha r_0 (= \alpha) &\\
     &  &  & \\
H^+ \rtimes \ZZ/2\ZZ & \overset{\sigma}{\longrightarrow} & W & \\
            r_i  & \longmapsto     & s_0 s_i & \text{ for }i=0,1,\ldots,n\\
          \alpha & \longmapsto     & s_0 &
\end{matrix}
$$
To check that $\rho$ is well-defined one must check that the $(W,S)$ Coxeter
relations $(s_i s_j)^{m_{ij}}=e$ for $0 \leq i \leq j \leq n$ map under $\rho$ to relations in
$H^+ \rtimes \ZZ/2\ZZ$.  Bearing in mind that $r_0=e$, this is checked as follows:
$$
(s_i s_j)^{m_{ij}}=e \quad \mapsto \quad \left( \alpha r_i \alpha r_j \right)^{m_{ij}}
                          = \left( r_i^{-1} \alpha \alpha r_j \right)^{m_{ij}} \\
                          = \left( r_i^{-1} r_j \right)^{m_{ij}} = e.
$$
To check that $\sigma$ is well-defined one can check that the
relations in the second presentation for $H^+ \rtimes \ZZ/2\ZZ$ map under $\sigma$ to
relations in $W$. These are checked as follows:
$$
\begin{aligned}
(r_i^{-1} r_j)^{m_{ij}}=e &\quad \mapsto \quad
    \left( s_i s_0 s_0 s_j \right)^{m_{ij}} = \left( s_i s_j \right)^{m_{ij}} = e.\\
\alpha^2 =e &\quad \mapsto \quad s_0^2 = e \\
\alpha r_i \alpha = r_i^{-1} & \quad \mapsto \quad s_0 (s_0 s_i ) s_0 = (s_0
s_i)^{-1}.
\end{aligned}
$$

\noindent
Once one knows that $\rho, \sigma$ are well-defined, it is easily checked that they
are inverse isomorphisms by checking this on generators.

Since $\sigma(H^+) \subseteq W^+$, and both $\sigma(H^+), W^+$
are subgroups of $W$ of index $2$, it must be that $\sigma(H^+)=W^+$.  Hence $\sigma$
restricts to the desired isomorphism between the abstractly presented group $H^+$ and $W^+$.
\end{proof}

\subsection{Length with respect to $R \cup R^{-1}$}
\label{general-length-subsection}

The maps $\rho, \sigma$ which appear in the proof of Proposition~\ref{Bourbaki-exercise}
lead to a nice interpretation for the length function of $W^+$ with respect to the symmetrized generating
set $R \cup R^{-1}$.

\begin{defn} \rm \ \\
Given a group $G$ and subset $A \subset G$, let $A^*$ denote the set
of all words ${\bf a}=(a_1,\ldots,a_\ell)$ with letters $a_i$ in $A$.
Let $A^{-1}:=\{a^{-1}: a \in A\}$.

Let $\ell_A(\cdot)$ denote the length function on $G$ with respect to the set $A$, that is,
$$
\ell_A(g):=\min\{\ell: g = a_1 a_2 \cdots a_\ell \text{ for some }a_i \in A\},
$$
where by convention, we set $\ell_A(g)=\infty$ if there are no such expressions for $g$.

Given an $A^*$-word ${\bf a}$ that factors $g$ in $G$, say that ${\bf a}$ is a
{\it reduced word} for $g$ if it achieves the minimum possible length $\ell_A(g)$.
\end{defn}

\begin{defn} \rm \ \\
\label{nu-definition}
Given a Coxeter system $(W,S)$ with $S=\{s_0,s_1,\ldots,s_n\}$ as before, let
$\nu(w)$ denote the minimum number of generators $s_j \neq s_0$ occurring in any
expression ${\bf s}=(s_{i_1}, \cdots, s_{i_\ell}) \in S^*$ that factors $w$ in $W$,
i.e. $w=s_{i_1} \cdots s_{i_\ell}$.
\end{defn}

\begin{proposition}
\label{symmetrized-generator-length}
For a Coxeter system $(W,S)$ with $S=\{s_0,s_1,\ldots,s_n\}$ as before, and
the presentation $(W^+,R)$ in \eqref{general-W-plus-presentation}, one
has
$$
\ell_{R \cup R^{-1}}(w)=\nu(w)
$$
for all $w \in W^+$.
\end{proposition}
\begin{proof}
Assume $w \in W^+$.  First we prove the inequality $\ell_{R \cup R^{-1}}(w) \geq \nu(w)$.
Given an $(R\cup R^{-1})^*$-word ${\bf r}$ that factors $w$ of the shortest
possible length $\ell_{R \cup R^{-1}}(w)$, apply
the map $\sigma$ from before
$$
\begin{aligned}
r_i &\mapsto s_0 s_i\\
r^{-1}_i &\mapsto s_i s_0
\end{aligned}
$$
to each letter and concatenate.
This gives an $S^*$-word {\bf s} that factors $w$,
having $\ell_{R \cup R^{-1}}(w)$ occurrences
of generators $s_j \neq s_0$.  Hence the minimum possible such number $\nu(w)$ must
be at most $\ell_{R \cup R^{-1}}(w)$.

Similarly we prove the opposite inequality $\ell_{R \cup R^{-1}}(w) \leq \nu(w)$.
Given an $S^*$-word ${\bf s}$ that factors $w$ with the minimum
number $\nu(w)$ of occurrences of generators $s_j \neq s_0$,
apply the map $\rho$ from before
$$
\begin{aligned}
s_i &\mapsto \alpha r_i \text{ for }i=1,\ldots,n\\
s_0 &\mapsto \alpha
\end{aligned}
$$
to each letter and concatenate.
This gives an $(R \cup \{\alpha\})^*$-word {\bf r} that factors $w$,
having $\nu(w)$ occurrences of generators $r_i$, and an {\it even} number of
occurrences of $\alpha$ (because $w \in W^+$ implies ${\bf s}$ has even length).
Repeatedly using the relation $\alpha r_i \alpha^{-1} = r_i^{-1}$, one can bring all
these evenly many occurrences of $\alpha$ in ${\bf r}$ to the right end of the word,
where they will cancel out
because $\alpha^2=1$.  This leaves an $(R \cup R^{-1})^*$ word factoring $w$, having length
$\nu(w)$.  Hence $\ell_{R \cup R^{-1}}(w) \leq \nu(w)$.
\end{proof}

For any $w \in W^+$, the proof of the inequality $\ell_{R \cup
R^{-1}}(w) \leq \nu(w)$ describes in two steps a map (which we
will also call $\rho$) from $S^*$-words ${\bf s}$ factoring $w$ to
$(R \cup R^{-1})^*$-words ${\bf r}$ factoring $w$. For future use,
we point out that this map has the following simple explicit
description :
\begin{enumerate}
\item[$\bullet$]
replace $s_i$ with $r_i$ for $i=2,3,\ldots,n$,
\item[$\bullet$]
replace $s_1$ with $r_1$ or $r_1^{-1}$, respectively,
depending upon whether the letter $s_1$ occurs in an even or odd position of ${\bf s}$, respectively, and
\item[$\bullet$] remove all occurrences of $s_0$.
\end{enumerate}

As an example,
$$
\begin{matrix}
\text{position: }&(1,  &2   &3   &4  &5   &6   &7    &8   &9       )\\
S^*-\text{word: }&(s_0,&s_2,&s_0,&s_1,&s_2,&s_0,& s_0,&s_3,&s_1,     \\
(R \cup \{\alpha\})^*-\text{word: }
  &(\alpha,&\alpha r_2,&\alpha,&\alpha r_1,&\alpha r_2,&\alpha,& \alpha,&\alpha r_3,&\alpha r_1) \\
(R \cup R^{-1})^*-\text{word: }&(    &r_2,&    &r_1,&r_2^{-1},&    &     &r_3,&r_1^{-1})
\end{matrix}
$$

\begin{proposition}
\label{explicit-lift-map}
The map just described coincides with the map
from $S^*$-words factoring $w$ to $(R \cup R^{-1})^*$-words factoring $w$
described in the proof of Proposition~\ref{symmetrized-generator-length}.
\end{proposition}
\begin{proof}
Note that an occurrence of $r_i$ in ${\bf r}$ which came from an occurrence of
$s_i$ in the $k^{th}$ position of ${\bf s}$ will
start with $k$ occurrences of $\alpha$ to its left in the $(R \cup \{\alpha\})^*$-word,
and each of these $\alpha$'s ``toggles'' it between
$r_i \leftrightarrow r_i^{-1}$ as that $\alpha$ moves past it to the right.
\end{proof}

\begin{example}\label{ltrmin=length} \rm \ \\
Let $(W,S)$ be the symmetric group $W=\sym_n$, with
$S=\{s_0,s_1,\ldots,s_{n-2}\}$ in which $s_i$ is the adjacent transposition $(i+1,i+2)$,
so $s_0=(1,2)$; this is the usual Coxeter system of type $A_{n-1}$.
Then the length in $W^+=\alt_n$ with respect to generating set
$R \cup R^{-1} = R \cup \{r_1^{-1}\}$ was considered in \cite{RegevRoichman},
where it was given the following explicit interpretation, reproven here for
the sake of completeness.

Given a permutation $w \in \sym_n$, let $\ltrmin(w)$ denote its
number of {\it left-to-right minima}, that is, the number of
$j\in\{2,3,\ldots,n\}$ satisfying $w(i) > w(j)$ for $1 \leq i <
j$. Let $\inv(w)$ denote its number of {\it inversions}, that is,
the number of pairs $(i,j)$ with $1 \leq i < j  \leq n$ and $w(i)
> w(j)$.  It is well-known \cite[Proposition 1.5.2]{BB} that the
Coxeter group length $\ell_S$ has the interpretation
$\ell_S(w)=\inv(w)$.

\begin{proposition}
\label{ltrmin--length}
For any $w \in \sym_n$, the maximum number of occurrences of $s_0$ in a reduced $S^*$-word for $w$
is  $\ltrmin(w)$.
Consequently,
$$
\begin{aligned}
\ell_{R \cup R^{-1}}(w) & =\ell_S(w)-\ltrmin(w)\\
                            & = \inv(w)-\ltrmin(w).
\end{aligned}
$$
\end{proposition}
\begin{proof}
For the first assertion, consider a reduced word ${\bf s}$ factoring $w$ as
sorting $w$ to the identity permutation $e$ by a sequence of adjacent transpositions.
During the process $\ltrmin$ can only go weakly downward, never up, and each time
one performs $s_0$, $\ltrmin$ goes down by one.  Since $\ltrmin(e)=0$, this
implies $\ltrmin(w)$ provides an upper bound on the number of occurrences of $s_0$ in ${\bf s}$.
On the other hand, one can produce such a sorting sequence for $w$ having exactly $\ltrmin(w)$
occurrences of $s_0$ as follows:  first move the letter $n$ step-by-step to the $n^{th}$ position,
then move the letter $n-1$ to the $(n-1)^{th}$ position, etc.  It's not hard to see that this
will use an $s_0$ exactly $\ltrmin(w)$ times.

For the second assertion, note by Proposition~\ref{symmetrized-generator-length}
that $\ell_{R \cup R^{-1}}(w)$ is the minimum number of $s_j \neq s_0$ in an
$S^*$-word factoring $w$.  However, by the {\it deletion condition} or Tits' solution
to the word problem for $(W,S)$, this minimum will be achieved by some {\it reduced} $S^*$-word
that factors $w$ (there exists such a reduced factorization for $w$ which is a subword
of the original factorization).  A reduced word achieving this minimum will have
exactly $\ltrmin(w)$ occurrences of $s_0$ by the first assertion, and will have $\ell_S(w)$ letters total,
so it will have $\ell_S(w)-\ltrmin(w)$ occurrences of $s_j \neq s_0$.
\end{proof}

In \cite{RegevRoichman} it was shown that for $(W,S)$ of type $A_{n-1}$ with $s_0$
a leaf node as above, one has
\begin{equation}
\label{typeA-length-gf}
\sum_{w \in W^+} q^{\ell_{R \cup R^{-1}}(w)}
=(1+2q)(1+q+2q^2)\cdots (1+q+q^2+\cdots+q^{n-3}+2q^{n-2}),
\end{equation}
and there are refinements of \eqref{typeA-length-gf} that
incorporate other statistics; see \cite[Proposition 5.7(2),
5.11(2)]{RegevRoichman}. The results of the current paper do not
recover this, and are in a sense, complementary-- they say more
about the case where $s_0$ is evenly-laced.
\end{example}

\subsection{Parabolic subgroup structure for $(W^+,R)$}
\label{general-parabolic-subsection}

The presentation \eqref{general-W-plus-presentation}
for $W^+$ with respect to the generating set $R$
likens $(W^+,R)$ to a Coxeter system, and suggests the following
definition.

\begin{defn} \rm \ \\
For any $J \subset R = \{r_1,\ldots,r_n\}$, the subgroup
$W^+_J=\langle J \rangle$ generated by $J$ inside $W^+$ will be
called a {\it (standard) parabolic subgroup}.
\end{defn}

The structure of parabolic subgroups $W_J$ for $(W,S)$ is an important part of the theory.
For $(W^+,R)$ one finds that its parabolic subgroups are closely tied to the parabolic
subgroups $W_J$ containing $s_0$, via the following map.

\begin{defn} \rm \ \\
\label{tau-definition}
Define $\tau: W \rightarrow W^+$ by
$$
\tau(w) :=
\begin{cases}
   w    & \text{ if } w \in W^+ \\
   ws_0 & \text{ if } w \not\in W^+
\end{cases}
$$
In other words, $\tau(w)$ is the unique element in the coset $wW_{\{s_0\}}=\{w,ws_0\}$
that lies in $W^+$.
\end{defn}

The following key property of $\tau$ is immediate from its definition.
\begin{proposition}
\label{tau-is-equivariant}
The set map $\tau: W \rightarrow W^+$ is equivariant for the
$W^+$-actions on $W, W^+$ by left-multiplication.
\end{proposition}

In fact, $\tau$ induces a $W^+$-equivariant bijection $W/W_{\{s_0\}} \rightarrow W^+$,
but we'll soon see that more is true.  Given any $J \subseteq S$ with $s_0 \in J$, let
$$
\tau(J):=\{r_i: s_0 \neq s_i \in J\}.
$$
Note that the map $J \mapsto \tau(J)$ is
a bijection between the indexing sets for parabolic subgroups in $W$ containing $s_0$
and for all parabolic subgroups of $W^+$.

\begin{proposition}
For any $J \subseteq S$ with $s_0 \in J$, one has $$
W_J \cap W^+ = W^+_{\tau(J)}.
$$
\end{proposition}
\begin{proof}
The inclusion  $W^+_{\tau(J)} \subseteq W_J \cap W^+$ should be clear.
For the reverse inclusion, given $w \in W_J \cap W^+$, write a $J^*$-word ${\bf s}$
that factors $w$ containing only $s_i \in J$.  Applying the map $\rho$ from
Proposition~\ref{explicit-lift-map} gives a $(\tau(J) \cup \tau(J)^{-1})^*$-word that
factors $w$, showing that $w \in W^+_{\tau(J)}$.
\end{proof}

\begin{proposition}
\label{tau-bijective}
For any $J \subseteq S$ with $s_0 \in J$, the (set) map $\tau$ induces a $W^+$-equivariant
bijection
$$
W/W_J \overset{\tau}{\longrightarrow} W^+/W^+_{\tau(J)}.
$$
In particular, taking $J=\{s_0\}$, this is a $W^+$-equivariant bijection
$$
W/W_{\{s_0\}} \overset{\tau}{\longrightarrow} W^+.
$$
\end{proposition}
\begin{proof}
One has a well-defined composite map of sets
$$
W \overset{\tau}{\rightarrow} W^+ \rightarrow W^+/W^+_{\tau(J)}
$$
sending $w$ to $\tau(w)W^+_{\tau(J)}$.  This composite surjects
because $\tau: W \rightarrow W^+$ surjects.

It remains to show two things:  the  composite induces a
{\it well-defined} map $W/W_J \overset{\tau}{\rightarrow} W^+/W^+_{\tau(J)}$,
and that this induced map is {\it injective}.  Both of these are shown
simultaneously as follows: for any $u,v \in W$ one has
\begin{eqnarray*}
\tau(u)W^+_{\tau(J)} = \tau(v)W^+_{\tau(J)}
  &\Leftrightarrow& \tau(v)^{-1}\tau(u) \in W^+_{\tau(J)} \\
  &\Leftrightarrow& \tau ( \tau(v)^{-1} u ) \in W^+_{\tau(J)} \\
  &\Leftrightarrow& \tau(v)^{-1} u \in W_J \\
  &\Leftrightarrow& v^{-1} u \in W_J \\
  &\Leftrightarrow& uW_J = vW_J
\end{eqnarray*}
where we have used throughout the fact that $s_0 \in J$, and where
the second equivalence uses the $W^+$-equivariance
of the set map $\tau: W \rightarrow W^+$ from Proposition~\ref{tau-is-equivariant}.
\end{proof}

Note that Proposition~\ref{tau-bijective} implies
that for any $J \subseteq S$ with $s_0 \in J$, the set of minimum $\ell_S$-length coset
representatives $W^J$ for $W/W_J$ maps under $\tau$ to a set
$\tau(W^J)$ of coset representatives for $W^+/W^+_{\tau(J)}$.  It turns out that
these coset representatives $\tau(W^J)$ are always of minimum $\ell_{R \cup R^{-1}}$-length.
To prove this, we note a simple property of the function $\nu$ that was defined
in Definition~\ref{nu-definition}.

\begin{proposition}
\label{nu-is-s0-invariant}
For any $w$ in $W$ one has
$$
\nu(s_0 w) = \nu(w)=\nu(ws_0)=\ell_{R \cup R^{-1}}(\tau(w)).
$$
\end{proposition}
\begin{proof}
Since $\nu(w^{-1}) = \nu(w)$, the first equality follows if one shows
the middle equality.  Also, since $\ell_{R \cup R^{-1}}(\tau(w)) = \nu(\tau(w))$ and
since $\tau(w)$ is either $w$ or $ws_0$, the last equality also follows
from the middle equality.

To prove the middle equality, it suffices to show the
inequality $\nu(ws_0) \leq \nu(w)$ for all $w \in W$;
the reverse inequality follows since $w= ws_0 \cdot s_0$.
But this inequality is clear: starting with an $S^*$-word ${\bf s}$ for $w$ that
has the minimum number $\nu(w)$ of occurrences of $s_j \neq s_0$, one can append an $s_0$ to
the end to get an $S^*$-word that factors $ws_0$ having no more such occurrences.
\end{proof}

\begin{corollary}
\label{achieve-min-length}
For any $J \subseteq S$ with $s_0 \in J$, the coset representatives $\tau(W^J)$ for $W^+/W^+_{\tau(J)}$
each achieve the minimum $\ell_{R \cup R^{-1}}$-length within their coset.
\end{corollary}
\begin{proof}
Let $w \in W^J$, and $w' \in \tau(w) W^+_{\tau(J)}$. Given an
$S^*$-word for $w'$ that has the minimum number $\nu(w')$ of
occurrences of $s_j \neq s_0$, one can extract from it an
$S^*$-reduced subword for $w'$.  Since $w'\in wW_J$ and $w \in
W^J$, one has $w \leq w'$ in the {\it strong Bruhat} order on $W$,
and hence one can extract from this a further $S^*$-subword
factoring $w$ \cite[\S 5.10]{Humphreys}. Consequently $\nu(w) \leq
\nu(w')$. But then Proposition~\ref{nu-is-s0-invariant} says that
$$
\ell_{R \cup R^{-1}}(\tau(w)) = \nu(w) \leq \nu(w') = \ell_{R \cup R^{-1}}(w')
$$
as desired.
\end{proof}

Note that we have made no assertion here about an element of
$\tau(W^J)$ being {\it unique} in achieving the minimum length
$\ell_{R \cup R^{-1}}$ within its coset, nor have we asserted that
the unique factorization $W^+ = \tau(W^J) \cdot W^+_{\tau(J)}$ has
additivity of lengths $\ell_{R \cup R^{-1}}$. In fact,
these properties fail in general
(see Remark~\ref{weak-descent-example}), but they will be shown in
Subsection~\ref{parabolic-cosets-revisited} to hold whenever $s_0$
is an evenly-laced node.

\begin{remark} \rm \ \\
The proof of Corollary~\ref{achieve-min-length} contains a fact which we
isolate here for future use.

\begin{proposition}
\label{Bruhat-nu-relation} Let $(W,S)$ be an arbitrary Coxeter
system. If $w < w'$ in the strong Bruhat order on $W$ then $\nu(w)
\le \nu(w')$. In particular,
\begin{enumerate}
\item[(i)] for any $s \in S, w \in W$, if $\ell_S(ws) < \ell_S(w)$
then $\nu(ws) \le \nu(w)$. \item[(ii)] for $w, w' \in W^+$, if $w
< w'$ in the strong Bruhat order on $W$ then $\ell_{R \cup
R^{-1}}(w) \le \ell_{R \cup R^{-1}}(w')$.
\end{enumerate}
\end{proposition}

\end{remark}

\subsection{The Coxeter complex for $(W^+,R)$}
\label{Coxeter-complex-subsection}

  Associated to every Coxeter system $(W,S)$ is a simplicial complex
$\Delta(W,S)$ known as its {\it Coxeter complex}, that has many
guises (see \cite[\S 1.15, 5.13]{Humphreys} and \cite[Exercise 3.16]{BB}):
\begin{enumerate}
\item[(i)] It is the {\it nerve} of the covering of the set $W$ by
the sets
$$
\{ wW_{S\setminus \{s\}} \}_{w \in W, s \in S}
$$
which are all cosets of maximal (proper) parabolic subgroups.
\item[(ii)] It is the unique simplicial complex whose face poset
has elements indexed by the collection
$$
\{ wW_J \}_{w \in W, J \subseteq S}
$$
of all cosets of all parabolic subgroups, with ordering by {\it reverse inclusion}.
\item[(iii)] It describes the decomposition by reflecting hyperplanes into cells
(actually spherical simplices) of the unit sphere intersected with the Tits cone in
the {\it contregredient representation} $V^*$ of $W$.
\end{enumerate}

The Coxeter complex $\Delta(W,S)$ enjoys many nice combinatorial, topological,
and representation-theoretic properties (see \cite{Bjorner}, \cite[Exercise 3.16]{BB}), such as:
\begin{enumerate}
\item[(i)] It is a pure $(|S|-1)$-dimensional simplicial complex, and is
{\it balanced} in the sense that if one colors the typical vertex
$\{ wW_{S\setminus \{s\}} \}_{w \in W, s \in S}$ by the element $s \in S$, then
every maximal face of $\Delta(W,S)$ contains exactly one vertex of
each color $s \in S$.

\item[(ii)] It is a {\it shellable pseudomanifold}, homeomomorphic to either an
$(|S|-1)$-dimensional sphere or open ball, depending upon whether $W$ is finite or infinite.

\item[(iii)]  When $W$ is finite, the homology $H_*(\Delta(W,S),\ZZ)$, which
is concentrated in the top dimension $|S|-1$, carries the {\it sign} character
of $W$.

\item[(iv)]  For each $J \subseteq S$, the
{\it type-selected} subcomplex $\Delta(W,S)_J$, induced on the
subset of vertices with colors in $J$, inherits the properties of being
pure $(|J|-1)$-dimensional, balanced, and shellable.  Consequently,
although $\Delta(W,S)_J$ is no longer homeomorphic
to a sphere, it is homotopy equivalent to a wedge of $(|J|-1)$-dimensional spheres.
Furthermore, the $W$-action on its top homology has an explicit decomposition into
Kazhdan-Lusztig cell representations.
\end{enumerate}

  The results of Section~\ref{general-parabolic-subsection}
allow us to define a Coxeter-like complex for $(W^+,R)$
in the sense of \cite{BabsonReiner}, and the map
$\tau$ allows one to immediately carry over many of the properties of $\Delta(W,S)$.

\begin{defn} \rm \ \\
Given a Coxeter system $(W,S)$ with $S=\{s_0,s_1,\ldots,s_n\}$, and the
ensuing presentation \eqref{general-W-plus-presentation} for $W^+$ via the
generators $R=\{r_1,\ldots,r_n\}$, define the {\it Coxeter complex} to be
the simplicial complex $\Delta(W^+,R)$ which is the nerve of the covering of the set $W^+$ by the
maximal (proper) parabolic subgroups
$$
\{ wW^+_{R \setminus \{r\}} \}_{w \in W^+, r \in R}.
$$
\end{defn}

Proposition~\ref{tau-bijective} and the usual properties of the
Coxeter complex $\Delta(W,S)$ immediately imply the following.

\begin{proposition}
\label{Coxeter-complex-isomorphism}
The map $\tau: W \rightarrow W^+$ induces a $W^+$-equvariant simplicial isomorphism
$$
\Delta(W,S)_{S\setminus\{s_0\}} \cong \Delta(W^+,R)
$$
where $\Delta(W,S)_{S\setminus\{s_0\}}$ denotes the type-slected subcomplex
obtained by deleting all vertices of color $s_0$ from $\Delta(W,S)$.

Consequently $\Delta(W^+,R)$ is a pure $(n-1)$-dimensional shellable
simplicial complex, which is balanced with color set $R$.

Similarly for any
$J \subseteq R$, its type-selected subcomplex $\Delta(W^+,R)_J$
is $W^+$-equivariantly isomorphic to the type-selected subcomplex
$\Delta(W,S)_{\tau^{-1}(J)}$.
\end{proposition}

This has consequences for the homology of $\Delta(W^+,R)$.
Let $\ZZ[W/W_{S-\{s_0\}}]$ denote the permutation action of $W^+$ on
cosets of the maximal parabolic $W_{S-\{s_0\}}$.  In other words,
$$
\ZZ[W/W_{S-\{s_0\}}] = \mathrm{Res}^W_{W^+} \mathrm{Ind}^W_{W_{S-\{s_0\}}} \mathbf{1}.
$$
If $W$ is finite, denote by $\ZZ v$ the unique copy of the trivial representation contained
inside $\ZZ[W/W_{S-\{s_0\}}]$, spanned by the sum $v$ of all cosets $wW_{S-\{s_0\}}$.

\begin{corollary}
The reduced homology $\tilde{H}_*(\Delta(W^+,R),\ZZ)$ is concentrated in dimension $n-1$,
and carries the
$W^+$-representation which is the restriction from $W$ of the representation
on the top homology of $\Delta(W,S)_{S\setminus\{s_0\}}$.
More concretely,
$$
H_*(\Delta(W^+,R),\ZZ) \cong
\begin{cases}
\ZZ[W/W_{S-\{s_0\}}] & \text{ when }W\text{ is infinite,} \\
\ZZ[W/W_{S-\{s_0\}}]/ \ZZ v & \text{ when }W\text{ is finite.}
\end{cases}
$$
\end{corollary}
\begin{proof}
The first assertions follow from Proposition~\ref{Coxeter-complex-isomorphism}
and the fact that a pure shellable $d$-dimensional complex has reduced homology concentrated in
dimension $d$.

The more concrete description of the $W^+$-action is derived as follows.
One can always apply {\it Alexander duality} to the embedding of
$\Delta(W,S)_{S-\{s_0\}}$ inside a certain $(|S|-1)$-dimensional sphere $\SSS^{|S|-1}$;
this sphere $\SSS^{|S|-1}$ is either
$\Delta(W,S)$ or its one-point compactification, depending upon whether $W$ is finite
or infinite.  In both cases, $W$ acts on the top homology $\tilde{H}_{|S|-1}(\SSS^{|S|-1},\ZZ)=\ZZ$
of this sphere by the sign character $\epsilon$, giving the following isomorphism of $W$-representations
(cf. \cite[Theorem 2.4]{Stanley}):
$$
\tilde{H}_{|S \setminus J|-1}(\Delta(W,S)_{S \setminus J},\ZZ)
\cong
\epsilon \otimes \left( \tilde{H}_{|J|-1}(\Delta(W,S)_{J},\ZZ) \right)^*.
$$
for any $J \subseteq S$; here $U^*$ denotes the contragredient
of a representation $U$, and when $W$ is infinite, the space $\Delta(W,S)_{J}$
appearing on the right should be replaced by its disjoint union $\Delta(W,S)_{J} \cup \{*\}$
with the compactification point $*$ of the sphere.

Taking $J=\{s_0\}$, one obtains a $W$-representation isomorphism between
the homology
$\tilde{H}_{|S|-2}(\Delta(W,S)_{S-\{s_0\}},\ZZ)$ and the twist by $\epsilon$
of either $\ZZ[W/W_{S-\{s_0\}}]$ or $\ZZ[W/W_{S-\{s_0\}}]/\ZZ v$, depending upon
whether $W$ is infinite or finite.   Restricting this isomorphism to $W^+$,
the twist by $\epsilon$ becomes trivial, and one gets the statement of
the corollary.
\end{proof}

\begin{example} \rm \ \\
Let $(W,S)$ be of type $A_3$, so that $W=\sym_3$, having Coxeter diagram which is
a path with three nodes.  If one labels the generators $S$ as
$$
\begin{aligned}
S&=\{s_0,s_1,s_2\}\\
&=\{(1,2),(2,3),(3,4)\},
\end{aligned}
$$
so that $s_0$ is a leaf node in the Coxeter diagram, then
Figure~\ref{Coxeter-complex-figure}(a) shows the Coxeter complex $\Delta(W^+,R)$ with facets
labelled by $W^+$.   Figure~\ref{Coxeter-complex-figure}(b) shows the isomorphic
type-selected subcomplex $\Delta(W,S)_{S-\{s_0\}}$ with facets labelled by $W^{\{s_0\}}$.

Figure~\ref{Coxeter-complex-figure}(c) shows the resulting Coxeter complex $\Delta(W^+,R)$ with facets
labelled by $W^+$ after one relabels
$$
\begin{aligned}
S&=\{s_0,s_1,s_2\}\\
&=\{(2,3),(1,2),(3,4)\},
\end{aligned}
$$
so that now $s_0$ is the central node, not a leaf, and $s_1, s_2$ commute.

\begin{figure}
\epsfxsize=120mm
\epsfbox{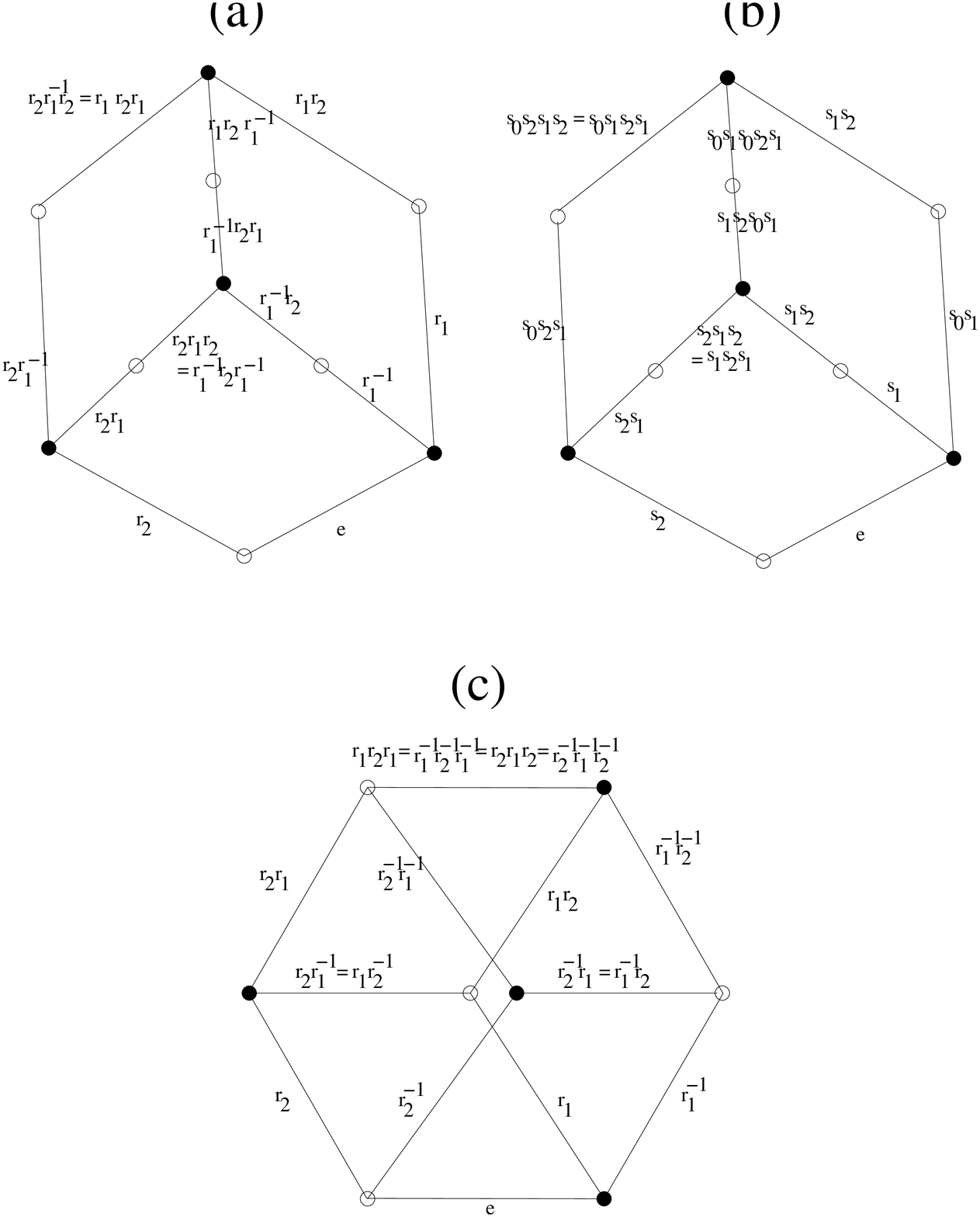}
\caption{Coxeter complexes for $(W^+,R)$ with $(W,S)$ of type
$A_3$, and two different choices for the distinguished node $s_0$.
Figure (a) shows $\Delta(W^+,R)$ when $s_0$ is a leaf node, that is, the
Coxeter diagram is labelled  $s_0 - s_1 - s_2$,
while (b) shows the isomorphic complex $\Delta(W,S)_{S-\{s_0\}}$. Figure (c) shows
$\Delta(W^+,R)$ when $s_0$ is the non-leaf node, that is, the
Coxeter diagram is labelled $s_1 - s_0 - s_2$.}
\label{Coxeter-complex-figure}
\end{figure}

\end{example}

\subsection{Palindromes versus reflections}
\label{palindrome-section}

  For a Coxeter system $(W,S)$, the set of reflections
$$
\T:=\bigcup_{\substack{w \in W \\ s \in S}} wsw^{-1}
$$
plays an important role in the theory.  A similar role for
$(W^+,R)$ is played by the set of {\it palindromes},  particularly
when $s_0$ is evenly-laced.  Palindromes will also give the
correct way to define the analogues of the strong Bruhat order
defined in Subsection~\ref{orders-subsection} below.

\begin{defn} \rm \ \\
Given a pair $(G,A)$ where $G$ is a group generated by a set $A$,
say that an element $g$ in $G$ is an {\it (odd) palindrome} if
there is an $(A \cup A^{-1})^*$-word ${\bf a}=(a_1,\ldots,a_\ell)$
factoring $g$ with $\ell$ odd such that $a_{\ell+1-i}=a_i$ for all $i$.
Denote the set of (odd) palindromes in $G$ by $\Pal(G)$.
\end{defn}

The set of palindromes for $(G,A)$ is always closed under taking inverses.
For a Coxeter system $(W,S)$, since $S$ consists entirely of involutions,
the set of palindromes is the same as the set $\T$ of reflections.

When $s_0$ is not evenly-laced in $(W,S)$, the palindromes
$\Pal(W^+)$ can behave unexpectedly, e.g.
the identity element $e$ is a palindrome: if $m_{01}$ is odd, one has the
odd palindromic expression $e = r_1^{m_{01}}$.
See also Example~\ref{odd-dihedral-example} below.

Nevertheless, one does have in general a very close relation between palindromes
for $(W^+,R)$ and palindromes (=reflections) for $(W,S)$.
Let
$$
\hat{\T} := \bigcup_{\substack{w \in W\\ s \in S \setminus \{s_0\}}} wsw^{-1}
$$
denote the set of reflections in $W$ that are conjugate to at least one
$s \neq s_0$.

\begin{proposition}
\label{reflection-subset}
The inclusion  $\hat{\T}\subset \T$ is proper if and
only if $s_0$ {\it is} evenly-laced.
\end{proposition}
\begin{proof}
When $s_0$ is not evenly laced, say $m_{01}$ is odd, then
$s_0$ is conjugate to $s_1$ and hence $\hat{\T} = \T$.

When $s_0$ {\it is} evenly-laced, the character $\chi_0: W \rightarrow \{\pm 1\}$ taking
value $-1$ on $s_0$ and $+1$ on $s_1,\ldots,s_n$ shows that
$s_0$ is not conjugate to any of $s_1,\ldots,s_n$, and hence the inclusion
$\hat{\T} \subsetneq \T$ is proper.
\end{proof}

\begin{proposition}
\label{general-palindrome-reflection-relation}
For any Coxeter system $(W,S)$, one has
$$
\Pal(W^+)s_0 = \hat{\T} = s_0\Pal(W^+)
$$
In other words, an element $w \in W^+$ is a palindrome with
respect to $R$ if and only if $ws_0$ (or equivalently $s_0w$) is a
reflection lying in the subset $\hat{\T}$, and vice-versa.
\end{proposition}
\begin{proof}
Since $\Pal(W^+)=\Pal(W^+)^{-1}$, it suffices to show the first equality.

Assume $w \in W^+$ is a palindrome, say
$w=r^{(1)} \cdots r^{(k-1)} r^{(k)} r^{(k-1)} \cdots r^{(1)}$ with each
$r^{(i)} \in R \cup R^{-1}$.  Then
$$
\begin{aligned}
ws_0 &=r^{(1)} \cdots r^{(k-1)} r^{(k)} r^{(k-1)} \cdots r^{(1)} s_0\\
     &=r^{(1)} \cdots r^{(k-1)} r^{(k)}s_0 (r^{(k-1)})^{-1} \cdots (r^{(1)})^{-1}\\
     &= u r^{(k)}s_0 u^{-1}
\end{aligned}
$$
for $u:=r^{(1)} \cdots r^{(k-1)}$, and where we have
used the fact that $r s_0 = s_0 r^{-1}$ for any
$r \in R \cup R^{-1}$.  Since $r^{(k)}s_0$ is either $s_0 s_i s_0$ or
$s_i s_0 s_0 = s_i$ for some $i=1,2,\ldots,n$, one concludes that $ws_0$
lies in $\hat{\T}$.

Conversely, given $ws_iw^{-1}$ in $\hat{\T}$, write any $S^*$-word ${\bf s}$ for
$w$.  Its reverse ${\bf s}^{rev}$ is a word for $w^{-1}$, and
$({\bf s}, s_i, {\bf s}^{rev},s_0)$ is a word for $ws_iw^{-1}s_0$.  Applying
the map from Proposition~\ref{explicit-lift-map} to this word yields an $(R \cup R^{-1})^*$
word ${\bf r}$ for $ws_iw^{-1}s_0$, which will be palindromic because there is
an odd distance in the word $({\bf s}, s_i, {\bf s}^{rev},s_0)$ between any two corresponding
occurrences of $s_j$ for $j=1,2,\ldots,n$.
\end{proof}

\begin{defn} \rm \ \\
Given $w \in W$, recall that its set of {\it left-shortening reflections} is
$$
\T_L(w):=\{ t \in \T: \ell_S(tw) < \ell_S(w) \}.
$$
Given $w \in W^+$, define its set of {\it left-shortening palindromes} by
$$
\Pal_L(w):=\{ p \in \Pal(W^+): \ell_{R \cup R^{-1}}(pw) < \ell_{R \cup R^{-1}}(w) \}.
$$
\end{defn}

In a Coxeter system $(W,S)$, it is well-known (\cite[Chapter 1]{BB}, \cite[\S 5.8]{Humphreys}
that for any $w$ in $W$, the set $\T_L(w)$ enjoys these properties:
\begin{enumerate}
\item[(a)] $\ell_S(w) = |\T_L(w)|$.

\item[(b)] ({\it strong exchange property}) For any $t \in \T$, and any
reduced $S^*$-word ${\bf s}=(s^{(1)},\ldots,s^{(\ell)})$ for $w$, the
following are equivalent
  \begin{enumerate}
  \item[(i)] $t \in \T_L(w)$, that is, $\ell_S(tw) < \ell_S(w)$.
  \item[(ii)] $t=t_k:=s^{(1)}\cdots  s^{(k-1)} s^{(k)} s^{(k-1)}\cdots s^{1)}$ for some $k$.
  \item[(iii)] $tw = s^{(1)}\cdots  s^{(k-1)} s^{(k+1)}\cdots s^{(\ell)}$ for some $k$.
  \end{enumerate}
 In other words, $\T_L(w) = \{t_1,\ldots,t_\ell\}$.

\item[(c)] The set $\T_L(w)$ determines $w$ uniquely.
\end{enumerate}

Analogously, given a reduced $(R \cup R^{-1})^*$-word ${\bf r} = (r^{(1)},\ldots,r^{(\nu(w))})$
that factors $w$ in $W^+$, one can define for $k=1,2,\ldots,\nu(w)$ the palindromes
$$
p_k:=(r^{(1)})^{-1} (r^{(2)})^{-1} \cdots (r^{(k)})^{-1} \cdots (r^{(2)})^{-1} (r^{(1)})^{-1}.
$$
One can relate this to $\Pal_L(w)$ and to $T_L(w)$ in general;  define for $w \in W$ the set
$$
\hat{\T}_L(w):=\T_L(w) \cap \hat{\T}.
$$

\begin{proposition}
\label{shortening-palindrome-reflection-relation}
For any choice of distinguished generator $s_0$, and for any $w \in W^+$, with the
above notation one has inclusions
\begin{equation}
\label{possibly-strict-inclusions}
\{ p_1,\ldots, p_{\nu(w)} \}
   \subseteq \Pal_L(w)
    \subseteq \hat{\T}_L(s_0w) s_0.
\end{equation}
When $s_0$ is evenly-laced, both inclusions are equalities:
$$
\{ p_1,\ldots, p_{\nu(w)} \} = \Pal_L(w) = \hat{\T}_L(s_0w) s_0.
$$

\end{proposition}
\begin{proof}
The first inclusion in \eqref{possibly-strict-inclusions} is straightforward, as
one calculates
$$
p_k w
 = (r^{(1)})^{-1} (r^{(2)})^{-1} \cdots (r^{(k-1)})^{-1} r^{(k+1)} r^{(k+2)} \cdots r^{(\nu(w))}
$$
and hence $\ell_{R \cup R^{-1}}(p_k w)< \nu(w) = \ell_{R \cup R^{-1}}(w)$.

For the second inclusion in \eqref{possibly-strict-inclusions},
given a palindrome $p\in \Pal(W^+)$, we know from Proposition~\ref{general-palindrome-reflection-relation}
that $t:=ps_0$ is a reflection in $\hat{\T}$, and conversely any reflection $t$ in $\hat{\T}$
will have $p:=ts_0$ a palindrome in $\Pal(W^+)$.  Thus it remains to show that
$$
\ell_{R \cup R^{-1}}(pw) < \ell_{R \cup R^{-1}}(w) \quad
\text{ implies } \quad
\ell_S(t s_0w) < \ell_S(s_0w).
$$
Using $\ell_{R \cup R^{-1}}=\nu$, along with the fact that
$\nu(s_0 w)=\nu(w)$ by Proposition~\ref{nu-is-s0-invariant}, and setting
$w':=s_0 w$, one can rewrite this desired implication as
\begin{equation}
\label{desired-palindrome-reflection-equivalence}
\nu(t w') < \nu(w') \quad \text{ implies }
\quad \ell_S(t w') < \ell_S(w').
\end{equation}
%
%
We show the contrapositive: if $\ell_S(t w') \geq \ell_S(w')$
then $t w'$ is greater than $w'$ in the Bruhat order on $W$, and
hence $\nu(tw') \geq \nu(w')$ by Proposition~\ref{Bruhat-nu-relation}.

For the assertions of equality, assuming $s_0$ is evenly-laced, it suffices to show that
the two sets $\{ p_1,\ldots, p_{\nu(w)} \}$ and
$\hat{\T}_L(s_0 w) s_0$ both have the same cardinality, namely $\nu(w)$.

For the first set, it suffices to show that $p_i \neq p_j$ for $1 \leq i < j \leq \nu(w)$.
Supposing $p_i=p_j$ for the sake of contradiction, one has
\begin{equation}
\label{rewriting-trick}
w=p_i^{-1} p_j w = r^{(1)} \cdots r^{(i-1)} (r^{(i+1)})^{-1} \cdots (r^{(j-1)})^{-1}
r^{(j+1)} \cdots r^{(\nu(w))}
\end{equation}
which gives the contradiction that $\ell_{R \cup R^{-1}}(w) < \nu(w)$.

For the second set, let $\ell:=\ell_S(s_0w)$
and choose a reduced $S^*$-word ${\bf s}=(s^{(1)},\ldots,s^{(\ell)})$ that
factors $s_0 w$.  Defining
$$
t_k=t^{(1)} t^{(2)} \cdots t^{(k)} \cdots t^{(2)} t^{(1)}
$$
for $1 \leq k \leq \ell$, one has \cite[Corollary 1.4.4]{BB},
\cite[\S 5.8]{Humphreys} that the $t_k$ are all distinct, and
$\T_L(s_0w):=\{ t_k \}_{1 \leq k \leq \nu(w)}$. Since $\nu(w)=\nu(s_0w)$, there will be exactly
$\nu(w)$ indices $\{i_1,\ldots,i_{\nu(w)}\}$ for which $s^{(i_j)}
\neq s_0$. As $t_k \in \hat{\T}$ if and only if $s^{(k)} \neq s_0$
(due to $s_0$ being evenly-laced), this means
$$
|\hat{\T}_L(s_0 w) s_0|=
|\hat{\T}_L(s_0 w)|=
|\T_L(s_0 w) \cap \hat{\T}| =
|\{ t_{i_1} , \ldots, t_{i_{\nu(w)}} \}|=
\nu(w).
$$
\end{proof}

\begin{example} \rm \ \\
\label{odd-dihedral-example}
When $(W,S)$ is the dihedral Coxeter system $I_2(m)$ in which $S=\{s_0,s_1\}$
with $m:=m_{01}$, then $(W^+,R)$ is simply the cyclic group of order $m$.
If one chooses $m$ to be odd, then every element $w \in W^+$ is a palindrome,
i.e. $\Pal(W^+)=W^+$, and one has
$$
\Pal(W^+)s_0=W^+s_0=\hat{\T}=\T.
$$
Furthermore, if one picks $m$ odd and sufficiently large, it
illustrates the potential bad behavior of palindromes when $s_0$ is not evenly-laced.  For example,
in this situation, $w = r_1^{-1} r_1^{-1}$ will have both inclusions strict in
\eqref{possibly-strict-inclusions}:
$$
\begin{matrix}
 \{ p_1,\ldots,p_{\nu(w)} \} & \subsetneq & \Pal_L(w) & \subsetneq & \hat{\T}_L(s_0w) s_0 \\
\Vert & &\Vert& &\Vert\\
 \{ r_1, \, r_1 r_1\}
   &  & \{r_1, \, r_1 r_1, \, r_1 r_1 r_1 \}
    &  & \{e, \, r_1, \, r_1 r_1, \, r_1 r_1 r_1, \, r_1 r_1 r_1 r_1  \}.
\end{matrix}
$$


This dihedral example also shows why replacing the set $\Pal(W^+)$
of palindromes for $(W^+,R)$ with the set of {\it conjugates of}
$R \cup R^{-1}$
\begin{equation}
\label{conjugates}
\bigcup_{\substack{w \in W^+\\ r \in R \cup R^{-1}}} w r w^{-1}
\end{equation}
would be the wrong thing to do:  in this example,
$W^+$ is cyclic and hence abelian, so that this set of conjugates
in \eqref{conjugates}
is no larger than $R \cup R^{-1}=\{r_1,r_1^{-1}\}$ itself!
\end{example}

Example~\ref{odd-dihedral-example} shows that the analogues for
palindromes in $(W^+,R)$ of the properties $\ell_S(w)=|\T_L(w)|$ and
the strong exchange property for reflections in $(W,S)$ can fail when $s_0$ is
not evenly-laced.  They do hold under the evenly-laced assumption--
see Theorem~\ref{strong-exchange-property} below,
which furthermore asserts that the set $\Pal_L(w)$ determines $w \in W^+$ uniquely
when $s_0$ is evenly-laced.  This raises the following question.

\begin{question}\ \\ 
When $s_0$ is chosen arbitrarily, does  $\Pal_L(w)$ determine $w
\in W^+$ uniquely?
\end{question}

\subsection{Weak and strong orders}
\label{orders-subsection}

  For a Coxeter system $(W,S)$ there are two related partial orders (the weak and strong
Bruhat orders) on $W$ which form graded posets with rank function $\ell_S$.  Here we define
analogues for $(W^+,R)$.

\begin{defn} \rm \ \\
Define the {\it (left) strong order} $\leq_{LS}$ on $W^+$ as the reflexive and
transitive closure of the relation $w \overset{p}{\rightarrow} pw$ if $p \in \Pal(W^+)$
and $\ell_{R \cup R^{-1}}(w) < \ell_{R \cup R^{-1}}(pw)$. Similarly define
the {\it (right) strong order} $\leq_{RS}$.

Define the {\it (left) weak order} $\leq_{LW}$ on $W^+$ as the reflexive and
transitive closure of the relation $w \lessdot_{LW} rw$ if $r \in R \cup R^{-1}$
and $\ell_{R \cup R^{-1}}(w)+1=\ell_{R \cup R^{-1}}(rw)$.
Similarly define the {\it (right) weak order} $\leq_{RW}$.
\end{defn}

Several things should be fairly clear from these definitions:
\begin{enumerate}
\item[(i)] Because these are reflexive transitive binary relations
on $W^+$ that are weaker than the partial ordering by the length
function $\ell_{R \cup R^{-1}}$, they are actually partial orders
on the set $W^+$. In other words, taking the transitive closure
creates no directed cycles. \item[(ii)] Because the map $w \mapsto
w^{-1}$ preserves the set of palindromes $\Pal(W^+)$ and the
length function $\ell_{R \cup R^{-1}}$, it induces an isomorphism
between the left and right versions of the two orders.
\item[(iii)] The identity $e \in W^+$ is the unique minimum
element in all of these orders. \item[(iv)] The (left, right,
resp.) strong order is stronger than the (left, right, resp.) weak
order. 
\item[(v)] For every $u, v\in W^+$,   $v \le_{RW} u$ implies
$P_L(v)\subseteq P_L(u)$.

\end{enumerate}

\begin{question} 
Does the inclusion $P_L(v)\subseteq P_L(u)$ imply $v \le_{RW} u$?
\end{question}

\begin{figure}
\epsfxsize=60mm
\epsfbox{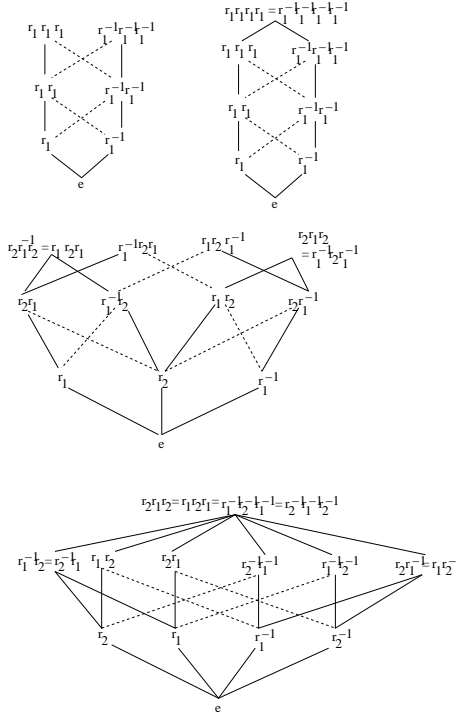}
\caption{Examples of the left weak (solid edges) and left strong orders (solid and
dotted edges) on $W^+$ for $(W,S)=I_2(7), I_2(8),$ and $A_3$ with $s_0$ labelling a leaf
node versus a non-leaf node.}
\label{Orderings-figure}
\end{figure}

Figure~\ref{Orderings-figure} shows the left weak and left strong orders
on $W^+$ for the two  dihedral Coxeter systems $I_2(7), I_2(8)$, as well as for
type $A_3$ with the two different choices for the node labelled $s_0$, as
in Figure~\ref{Coxeter-complex-figure}.

The usual weak and strong orders on a Coxeter group $W$ have several good properties
(see \cite[Chapters 2,3]{BB}):
they are all {\it graded} by the length function $\ell_S$, the left and right weak orders
are both {\it meet semilattices}, and the strong order is {\it shellable}.
A glance at Figure~\ref{Orderings-figure}
then raises several obvious questions about the analogous orders on $W^+$.

\begin{question} 
Are all of these orders graded by the function $\ell_{R \cup
R^{-1}}$, that is, do all maximal chains have the same length?
\end{question}

\begin{question} 
Do the weak orders form a meet semilattice in general?
\end{question}

\begin{question} 
Is the strong order shellable?
\end{question}

We will see in Subsection~\ref{orders-revisited}
that the answers to all of these questions are affirmative
when $s_0$ is evenly-laced.  Furthermore, in
Section~\ref{even-leaf-section} it will be shown that when $s_0$ is an evenly-laced
leaf node, the strong and weak orders coincide with the usual Coxeter group strong and weak orders
for the related Coxeter system $(W',S')$ defined there.

\begin{remark} \rm \ \\
Some things are clearly {\it not} true of the various orders, even in the best possible
situation where $s_0$ is an even leaf.

Although the left weak/strong orders are isomorphic to the right
weak/strong orders, they are not the {\it same} orders.  For
example, when $(W,S)$ is of type $B_3$ with $s_0$ the even leaf as
in Section~\ref{large-example} below, one can check that $r_1$ is
below $r_2 r_1$ in both the left weak and left strong orders, but
this fails in both the right weak and right strong orders.

None of the four orders (left/right weak/strong) on $W^+$
coincides with the restriction from $W$ to $W^+$ of the analogous
left/right weak/strong order on $W$.  For example, suppose that
$(W,S)$ has $W$ finite with an odd number $|T|$ of reflections,
and $s_0$ is evenly-laced (this occurs in type $B_n$ for $n$ odd; see
Section~\ref{large-example} below for the example of type $B_3$).  Then there will be
a maximum element, namely $\tau(w_0)=w_0s_0= s_0w_0$, for all four orders on $W^+$,
where here $w_0$ is the longest element in $W$; see Proposition~\ref{finite-case-maximum}
below. But $\tau(w_0)$ will {\it not} be a maximum element when one restricts any
of the left/right weak or strong orders from $W$ to $W^+$:  the elements $w_0 s_j$
for $j \geq 1$ will also lie in $W^+$, and have the same length
$$
\ell_S(w_0s_j)=|T|-1=\ell_S(\tau(w_0))
$$
and hence will be incomparable to $\tau(w_0)$.

Similarly, none of the four orders on $W^+$ coincides, via the bijection
$\tau: W^+ \rightarrow W^{\{s_0\}}$, to the restriction from $W$ to $W^{\{s_0\}}$
of the analogous order on $W$.  This can be seen already for $(W,S)$ of type $I_2(4)=B_2$,
where all four orders on $W^+$ are isomorphic to a rank two Boolean lattice,
while the various strong/weak orders restricted from $W$ to $W^{\{s_0\}}$
turn out either to be total orders or non-lattices.
\end{remark}

\section{The case of an evenly-laced node}
\label{evenly-laced-section}

  When the distinguished generator $s_0$ in $S=\{s_0,s_1,\ldots,s_n\}$
for the Coxeter system $(W,S)$ has the extra property that $m_{0i}$ is even
for $i=1,2,\ldots,n$, we say that $s_0$ is an {\it evenly-laced node}
of the Coxeter diagram.   This has many good consequences for the
presentation $(W^+,R)$ explored in the next few subsections:
\begin{enumerate}
\item[$\bullet$]
the length function $\ell_{R \cup R^{-1}}$ simplifies,
\item[$\bullet$]
the coset representatives $\tau(W^J)$ for $W^+/W_{\tau(J)}$ from Section~\ref{general-parabolic-subsection}
are distinguished by their minimum length within the coset, and the
length is additive in the decomposition $W^+ = \tau(W^J) \cdot W^+_{\tau(J)}$,
\item[$\bullet$]
the palindromes $\Pal(W^+)$ behave more like reflections, satisfying a
{\it strong exchange condition}, and consequently
\item[$\bullet$]
the partial orders considered earlier are as well-behaved as their analogues for
$(W,S)$.
\end{enumerate}

\subsection{Length revisited}
\label{evenly-laced-length-section}

\begin{defn} \rm \ \\
\label{nu-simplification-defn}
Part of Tits' solution to the word problem for the
Coxeter system $(W,S)$ asserts \cite[\S 3.3]{BB}
that one can connect any two reduced $S^*$-words for $w$ in $W$
by a sequence of {\it braid moves} of the form
\begin{equation}
\label{braid-move}
\underbrace{s_i s_j s_i s_j \cdots}_{m_{ij}\text{ letters}} =
\underbrace{s_j s_i s_j s_i \cdots}_{m_{ij}\text{ letters}}.
\end{equation}
When $s_0$ is evenly-laced, there will always be the same number of occurrences of $s_0$
on either side of \eqref{braid-move}, and hence the number of occurrences of $s_0$ in {\it any}
reduced word is the same;  denote this quantity $\ell_0(w)$.

The Coxeter presentation for $(W,S)$ also allows one to define, when $s_0$ is evenly-laced,
a homomorphism
\begin{equation}
\label{chi-zero-definition}
\begin{aligned}
\chi_0: W &\rightarrow \{\pm 1\}\\
\chi_0(s_0)&=-1 \\
\chi_0(s_j)&=+1 \text{ for }j=1,2,\ldots,n
\end{aligned}
\end{equation}
Note that $\chi_0(w) = (-1)^{\ell_0(w)}$.
\end{defn}

Recalling that $\nu(w)$ was defined to be the minimum number of $s_j \neq s_0$ occurring in an $S^*$-word that
factors $w$, one immediately concludes the following reinterpretation for the length function of $(W^+,R)$.

\begin{proposition}
\label{nu-simplification-prop}
Assume $s_0$ is evenly-laced.  Then for every $w \in W$ one has
$$
\nu(w)=\ell_S(w) - \ell_0(w).
$$
Consequently, for any $w \in W^+$, the length function
$\ell_{R \cup R^{-1}}(w) (=\nu(w))$ can be computed from any reduced $S^*$-word for $w$.
\end{proposition}

\subsection{Length generating function}
\label{length-gf-section}

  When $s_0$ is evenly-laced, the simpler interpretation for
the length function $\ell_{R \cup R^{-1}}$ allows one to compute its generating function
for $(W^+,R)$, by relating it to known variations on the usual Coxeter group
length generating function for $(W,S)$.

  The usual diagram-recursion methods \cite[\S 5.12]{Humphreys} for
writing down the Poincar\'e series
$$
W(S;q):=\sum_{w \in W} q^{\ell_S(w)}
$$
as a rational function in $q$ turn out to generalize
straightforwardly \cite{Macdonald, Reiner}, allowing one to write
down the finer Poincar\'e series
$$
W(S;q_0,q):=\sum_{w \in W} q_0^{\ell_0(w)} q^{\nu(w)}.
$$
This power series in $q_0,q$ will actually end up being a rational
function of $q_0,q$ for any Coxeter system $(W,S)$ with $s_0$ evenly-laced.
The key point is that in the unique factorization
$$
W = W^J  \cdot W_J,
$$
both statistics $\ell_0(w), \nu(w)$ behave additively
(see \cite[\S 5.12]{Humphreys} or \cite{Macdonald, Reiner}),
yielding the factorization
$$
W(S;q_0,q) = W^J(S;q_0,q)\cdot W_J(S;q_0,q).
$$
Here we are using the notation for any subset $A \subset W$ that
$$
A(S;q_0,q):=\sum_{w \in A}q_0^{\ell_0(w)} q^{\nu(w)}.
$$

\begin{defn} \rm \ \\
Define the $\ell_{R \cup R^{-1}}$ length generating function on $W^+$:
$$
W^+(R \cup R^{-1};q):=\sum_{ w \in W^+ } q^{\ell_{R \cup R^{-1}}(w)}.
$$
\end{defn}

\begin{corollary}\label{formula}
\label{evenly-laced-length-gf}
When $s_0$ is evenly-laced,
$$
W^+(R \cup R^{-1};q) = \left[ W^{\{s_0\}}(S;q_0,q) \right]_{q_0=1} =
\left[ \frac{W(S;q_0,q)}{1+q_0} \right]_{q_0=1}.
$$
\end{corollary}
\begin{proof}
Since the map $\tau: W^{\{s_0\}} \rightarrow W^+$ is a bijection
by Proposition~\ref{tau-bijective}, and since $\ell_{R \cup R^{-1}}(\tau(w)) = \nu(w)$
by Proposition~\ref{symmetrized-generator-length}, one has
$$
\begin{aligned}
W^+(R \cup R^{-1};q) &= \left[ W^{\{s_0\}}(S;q_0,q) \right]_{q_0=1}\\
         &= \left[ \frac{W(S;q_0,q)}{W_{\{s_0\}}(S;q_0,q)} \right]_{q_0=1}\\
         &= \left[ \frac{W(S;q_0,q)}{1+q_0} \right]_{q_0=1}.\\
\end{aligned}
$$
\end{proof}

\begin{example} \rm \ \\
\label{first-type-B-example}
Let $(W,S)$ be the Coxeter system of type $B_n(=C_n)$, so that $W$ is the group of
{\it signed permutations} acting on $\RR^n$.
Index $S=\{s_0,s_1,\ldots,s_{n-1}\}$ so that $s_0$ is the special generator that negates the
first coordinate, and $s_i$ swaps the $i^{th},(i+1)^{st}$ coordinates when $i\geq 1$.
The Coxeter presentation has
$$
\begin{aligned}
m_{01}&=4\\
m_{i,i+1}&=3 \text{ for }i=1,2,\ldots,n-1\\
m_{ij}&=2 \text{ for }|i-j| \geq 2.
\end{aligned}
$$
It is well-known
(see \cite{Bre, Macdonald, Reiner}) and not hard to check that
$$
W(S;q_0,q) = (-q_0;q)_n [n]!_q
$$
where
$$
\begin{aligned}
(x;q)_n&:=(1-x)(1-xq)(1-xq^2)\cdots(1-xq^{n-1})\\
[n]!_q &:=\frac{(q;q)_n}{(1-q)^n} = [n]_q [n-1]_q \cdots [2]_q [1]_q \\
[n]_q&:=\frac{1-q^n}{1-q} = 1+q+q^2 +\cdots +q^{n-1}.
\end{aligned}
$$
Consequently, Corollary~\ref{evenly-laced-length-gf} implies
$$
\begin{aligned}
W^+(R \cup R^{-1};q)
& =  \left[ \frac{(-q_0;q)_n [n]!_q}{1+q_0} \right]_{q_0=1} \\
& =  (-q;q)_{n-1} [n]!_q \\
& =[n]_q \prod_{j=1}^{n-1} (1+q^j) [j]_q.\\
& =[n]_q \prod_{j=1}^{n-1} [2j]_q.
\end{aligned}
$$
The same formula will be derived differently in Example~\ref{type-B-length-gf-again}.
\end{example}

\subsection{Parabolic coset representatives revisited}
\label{parabolic-cosets-revisited}

  Recall that for any subset $J \subseteq S$ with $s_0 \in J$, the
map $\tau$ sends the distinguished minimum $\ell_S$-length coset representatives
$W^J$ for $W/W_J$ to a collection $\tau(W^J)$ of coset representatives
for $W^+/W_{\tau(J)}$, each of which achieves the minimum $\ell_{R \cup R^{-1}}$-length
in its coset.  Thus for every $w \in W^+$ one has a unique factorization
\begin{equation}
\label{unique-factorization}
w = \tau(x) y
\end{equation}
with $x \in W^J$ and $y \in W^+_{\tau(J)}$ unique.  One can make a stronger
assertion when $s_0$ is evenly-laced.

\begin{proposition}
\label{parabolic-additivity}
Assume $s_0$ is evenly-laced.  Then in the unique factorization \eqref{unique-factorization}
one has additivity of lengths:
$$
\ell_{R \cup R^{-1}}(w) = \ell_{R \cup R^{-1}}(\tau(x)) + \ell_{R \cup R^{-1}}(y).
$$
\end{proposition}
\begin{proof}
Since elements $w \in W^+$ have $\ell_{R \cup R^{-1}}(w) = \nu(w)$, one must show that
in the factorization \eqref{unique-factorization}, one has
\begin{equation}
\label{desired-additivity}
\nu(w) = \nu(\tau(x)) + \nu(y).
\end{equation}

Because $s_0$ is evenly-laced, Definition~\ref{nu-simplification-defn}
and Proposition~\ref{nu-simplification-prop} imply that in the
usual length-additive parabolic factorization for $w \in W$ as
$w = w^J w_J$ with $w^J \in W^J, w_J \in W_J$, one has additivity of $\nu$:
\begin{equation}
\label{nu-additivity}
\nu(w) = \nu(w^J) + \nu(w_J).
\end{equation}

Note that \eqref{unique-factorization} implies that the usual parabolic factorization
$w=w^J \cdot w_J$ in $W$ must either take the form
$w = x \cdot y$ (if $\tau(x)=x$) or the form
$w= x \cdot s_0 y$ (if $\tau(x)= x s_0$).  In either case, the desired additivity \eqref{desired-additivity}
follows from \eqref{nu-additivity}, using Proposition~\ref{nu-is-s0-invariant}.
\end{proof}

This immediately implies the following.

\begin{corollary}\label{parabolic-unique}
When $s_0$ is evenly-laced, the coset representatives $\tau(W^J)$
for $W^+/W^+_{\tau(J)}$ can be distinguished intrinsically in any
of the following ways:
\begin{enumerate}
\item[(i)] $\tau(W^J)$ are the unique representatives within each
coset $wW^+_{\tau(J)}$ achieving the minimum $\ell_{R \cup R^{-1}}$-length.
\item[(ii)]
$$
\tau(W^J):=\{ x \in W^+: \ell_{R \cup R^{-1}}(xy) > \ell_{R \cup R^{-1}}(x)
\text{ for all }y \in W^+_{\tau(J)} \}.
$$
\item[(iii)]
$$
\tau(W^J):=\{ x \in W^+: \ell_{R \cup R^{-1}}(xr) > \ell_{R \cup R^{-1}}(x)
\text{ for all }r \in \tau(J) \cup \tau(J)^{-1} \}.
$$
\end{enumerate}
\end{corollary}

One also has the following immediate corollary, giving a factorization
for the $\ell_{R \cup R^{-1}}$ generating function.  Define the notation
for any subset $A \subset W^+$ that
$$
A(R\cup R^{-1};q):=\sum_{w \in A} q^{\ell_{R\cup R^{-1}}(w)}.
$$

\begin{corollary}\label{Poincare-factorization}
For every subset $J\subseteq R$
$$
W^+(R\cup R^{-1};q)= {W^+}^J(R\cup R^{-1};q) \cdot W^+_J(R\cup
R^{-1};q).
$$
\end{corollary}

Note that the factorization in Corollary~\ref{Poincare-factorization} fails
in general when $s_0$ is not evenly-laced.  For example, in the case of type $A_{n-1}$
where $W=\sym_n$ and $s_0$ is a leaf
node of the Coxeter diagram, $W^+(R\cup R^{-1};q)$ was given explicitly earlier in factored
form as \eqref{typeA-length-gf}, but is not divisible by $W^+_{\{r_i\}}(R \cup R^{-1};q)=1+q$
for any of the generators $r_i$ with $i > 1$.  See also Example~\ref{weak-descent-example} below.

\subsection{Descent statistics}
\label{first-descents-section}

  For a Coxeter system $(W,S)$, aside from the length statistic $\ell_S(w)$ for $w \in W$,
one often considers the {\it descent set} and {\it descent number} of $w$ defined by
$$
\begin{aligned}
\Des_S(w) &:=\{s \in S: \ell_S(ws) < \ell_S(w) \} \, \subseteq \, S\\
\des_S(w) &:= |\Des_S(w)|.
\end{aligned}
$$
Generating functions counting $W$ jointly by $\ell_S$ and $\Des_S(w)$ are
discussed in \cite{Reiner}.

  When $(W,S)$ is arbitrary, for the alternating group $W^+$ and its generating set $R$
there are several reasonable versions of the descent set one might consider.

\begin{defn} \rm \ \\
Given $w \in W^+$, define its {\it descent set} $\Des_{R \cup
R^{-1}}(w)$, {\it symmetrized descent set} $\widehat\Des_R(w)$,
{\it weak descent set}
(or {\it nonascent set}) $\WDes_{R \cup
R^{-1}}(w)$ and its {\it symmetrized  weak descent set}
$\HatWDes_R(w)$ as follows:
$$
\begin{aligned}
\Des_{R \cup R^{-1}}(w)
  &:=\{r \in R \cup R^{-1}: \ell_{R \cup R^{-1}}(wr) <  \ell_{R \cup R^{-1}}(w)\}
    \subseteq R \cup R^{-1} \\
\widehat\Des_R(w)
  &:=\{r \in R: \text{ either }r\text{ or }r^{-1}\in \Des_{R \cup R^{-1}}(w)\}
     \subseteq R \\
\WDes_{R \cup R^{-1}}(w)
  &:=\{r \in R \cup R^{-1}: \ell_{R \cup R^{-1}}(wr) \leq  \ell_{R \cup R^{-1}}(w)\}
    \subseteq R \cup R^{-1} \\
\HatWDes_R(w)
  &:=\{r \in R: \text{ either }r\text{ or }r^{-1}\in \WDes_{R \cup R^{-1}}(w)\}
     \subseteq R \\
\end{aligned}
$$
\end{defn}

Part of the justification for considering weak descents
comes from the type $A_{n-1}$ example where $W=\sym_n$:  in \cite[Theorem 1.10(2)]{RegevRoichman},
it was shown that the resulting major index
(i.e., the sum of the indices of the weak descents) is
equi-distributed with the length $\ell_{R\cup R^{-1}}$.

Note that one did not have to worry about {\it weak} descents for
$(W,S)$ because the existence of the sign character
shows that one always has $\ell_S(ws) \neq \ell_S(w)$ for any $s \in S$.
This can fail for $(W^+,R)$ and the length function $\ell_{R \cup
R^{-1}}$ in general.

\begin{example} \rm \ \\
\label{weak-descent-example}
Continuing Example~\ref{odd-dihedral-example},
let $(W,S)$ be the dihedral Coxeter
system $I_2(m)$ with $m=2k+1$.  Then the two elements
$r_1^k, r_1^{-k}$ both achieve the maximum
$\ell_{R \cup R^{-1}}$-length value of $k$, but differ by multiplication on the right
by elements of $R \cup R^{-1}$:
$$
\begin{aligned}
r_1^k \cdot &r_1 = r_1^{-k} \\
r_1^{-k} \cdot &r^{-1}_1 = r_1^k.
\end{aligned}
$$
Note that this also illustrates the failure of both Proposition~\ref{parabolic-additivity}
and Corollary~\ref{parabolic-unique} without the assumption that $s_0$
is evenly-laced: they fail on the coset
$r_1^k W^+_{\tau(J)}=r_1^{-k} W^+_{\tau(J)}$,
where $J=\{s_0,s_1\}$ and $\tau(J)=\{r_1\}$.

\end{example}

When $s_0$ is an evenly-laced node, restricting the
character $\chi_0$ to $W^+$ one has
$$
\chi_0(w) = (-1)^{\nu(w)} = (-1)^{\ell_{R \cup R^{-1}}(w)}.
$$
This shows that $\ell_{R \cup R^{-1}}(wr) \neq \ell_{R \cup R^{-1}}(w)$ for any $r \in R$,
and hence, in this case, weak descents are the same as descents:
$$
\begin{aligned}
\WDes_{R \cup R^{-1}}(w) =\Des_{R \cup R^{-1}}(w)
    &= \{r \in R \cup R^{-1}: \ell_{R \cup R^{-1}}(wr) <  \ell_{R \cup R^{-1}}(w)\}\\
\HatWDes_R(w)=\widehat\Des_R(w)
  &=\{r \in R: \text{ either }r\text{ or }r^{-1}\in \Des_{R \cup R^{-1}}(w)\}
\end{aligned}
$$

Note also that the set $\WDes_{R \cup R^{-1}}(w)$ completely
determines the set $\HatWDes_R(w)$, and hence is finer information
about $w$.  It would be nice to have generating functions
counting $W^+$ jointly by $\ell_{R\cup R^{-1}}$ and either
$\WDes_{R\cup R^{-1}}$ or $\HatWDes_R$.  These seem hard to
produce in general.
However, when $s_0$ is
evenly-laced, we next show how to produce such a generating
function for the pair $(\ell_{R\cup R^{-1}}, \HatWDes_R)$. In
Subsection~\ref{even-leaf-section}, we will do the same for the
finer information $(\ell_{R\cup R^{-1}}, \WDes_{R \cup R^{-1}})$
under the stronger hypothesis that $s_0$ is an evenly-laced {\it
leaf}.

It turns out that nonascents in $(W^+,R)$ relate to descents in
$(W,S)$ of the minimum length parabolic coset representatives
$W^{\{s_0\}}$ for $W/W_{\{s_0\}}$. This is mediated by the inverse
$\tau^{-1}$ to the bijection $\tau: W^{\{s_0\}} \rightarrow W^+$
that comes from taking $J=\{s_0\}$ in
Proposition~\ref{tau-bijective}.

  Our starting point is a relation for general $(W,S)$ between
$\HatWDes_R$ on $W^+$ and $\Des_S$ on $W^{\{s_0\}}$.
For the purpose of comparing subsets of $R=\{r_1,\ldots,r_n\}$
and $S \setminus \{s_0\}=\{s_1,\ldots,s_n\}$, identify both of
these sets of generators with their subscripts $[n]:=\{1,2,\ldots,n\}$.

\begin{proposition}
\label{descent-inclusion-equality-prop}
After the above identification of subscripts,
for any Coxeter system $(W,S)$ and $s_0 \in S$ and $w \in W^+$, one
has a (possibly proper) inclusion
\begin{equation}\label{descents-inclusion}
\HatWDes_R(w) \supseteq \Des_S( \tau^{-1}(w)).
\end{equation}
When $s_0$ is evenly-laced, this inclusion becomes an equality:
\begin{equation}\label{descents-equality}
 (\Des_R(w)= )\,\, \HatWDes_R(w) = \Des_S( \tau^{-1}(w)).
\end{equation}
\end{proposition}
\begin{proof}
To show the inclusion, given $w \in W^+$, assume $s_j \in \Des_S(\tau^{-1}(w))$,
and one must show that $r_j \in \HatWDes_R(w)$ (note that $j \neq 0$ since $\tau^{-1}(w) \in W^{\{s_0\}}$).
Since $\ell_S( \tau^{-1}(w)s_j) < \ell_S(\tau^{-1}(w))$,
by Proposition~\ref{Bruhat-nu-relation}(i) one has
$$
\nu( \tau^{-1}(w)s_j ) \leq \nu(\tau^{-1}(w)).
$$
If $\tau^{-1}(w)=w$ then this gives
$$
\ell_{R \cup R^{-1}}(w r_j^{-1})
 = \nu(w s_j s_0)
 = \nu(w s_j)
 \leq \nu(w)
 = \ell_{R \cup R^{-1}}(w)
$$
using Proposition~\ref{nu-is-s0-invariant}.
If $\tau^{-1}(w)=ws_0$ then
$$
\ell_{R \cup R^{-1}}(w r_j)
 = \nu(w s_0 s_j)
 \leq \nu(w s_0)
 = \nu(w)
 = \ell_{R \cup R^{-1}}(w)
$$
again using Proposition~\ref{nu-is-s0-invariant}.
Either way, one has $r_j \in \HatWDes(w)$.

Now assume $s_0$ is evenly-laced, and $r_j \in \HatWDes_R(w)( =
\widehat\Des_{R}(w))$.  One must show that $s_j \in
\Des_S(\tau^{-1}(w))$. Consider these cases:

\noindent {\bf Case 1.} $r_j \in \Des_{R\cup R^{-1}}(w)$. Then
$$
\nu(ws_0s_j) = \ell_{R \cup R^{-1}}(wr_j) < \ell_{R \cup R^{-1}}(w) = \nu(w) = \nu(ws_0),
$$
which forces $\ell_S(ws_0s_j) < \ell_S(ws_0)$ by
Proposition~\ref{Bruhat-nu-relation}(i). Thus $s_j \in \Des_S(ws_0)$.

If $\tau^{-1}(w) = ws_0$, then we're done.  If $\tau^{-1}(w) =w$, and one
assumes for the sake of contradiction that $s_j \not\in \Des_S(w)$,
then one has
$$
w \in W^{\{s_0\}} \cap  W^{\{s_j\}} =  W^{\{s_0,s_j\}}.
$$
This gives the contradiction
$$
\ell_S(ws_0s_j) = \ell_S(w) + 2 \not < \ell_S(w) + 1 = \ell_S(ws_0).
$$

\noindent {\bf Case 2.} $r^{-1}_j \in \Des_{R\cup R^{-1}}(w)$.
Then
$$
\nu(ws_j) = \nu(ws_js_0) = \ell_{R \cup R^{-1}}(wr_j^{-1}) < \ell_{R \cup R^{-1}}(w) = \nu(w),
$$
which forces $\ell_S(ws_j) < \ell_S(w)$ by Proposition~\ref{Bruhat-nu-relation}(i).
Thus $s_j \in \Des_S(w)$.

If $\tau^{-1}(w) = w$, then we're done.  If $\tau^{-1}(w) =ws_0$, and one
assumes for the sake of contradiction that $s_j \not\in \Des_S(ws_0)$,
then one has
$$
ws_0 \in W^{\{s_0\}} \cap  W^{\{s_j\}} =  W^{\{s_0,s_j\}}.
$$
This gives the contradiction
$$
\ell_S(ws_j) = \ell_S(ws_0 \cdot s_0s_j) = \ell_S(ws_0) +
2=\ell_S(w)+1 \not < \ell_S(w).
$$
\end{proof}

\begin{remark} \rm \ \\
To see that the inclusion in \eqref{descents-inclusion}
can be proper, consider the Coxeter
system $(W,S)$ of type $A_3$ with $s_0$ chosen to be a leaf node,
as in Figure~\ref{Coxeter-complex-figure}(a). Here if one takes
$w=r^{-1} r_2 r_1$ then $\tau^{-1}(w)=s_1 s_2 s_0 s_1$, with
$\HatWDes_R(w)=\{r_1,r_2\}$ but $\Des_S(\tau^{-1}(w))=\{s_1\}$.

We should also point out that this problem cannot
be fixed by using $\widehat\Des_R$ instead of $\HatWDes_R(w)$.
Not only would this not give equality in \eqref{descents-inclusion},
but one would no longer in general have an inclusion:
For example, for $w=r_1 r_2 r_1^{-1}\in A_3$
with $s_0$ chosen to be a leaf node as above,
$$
\begin{aligned}
\widehat\Des_R(w)&=\{r_1\}\\
\Des_S (\tau^{-1}(w))&=\Des_S(s_0s_1s_0s_2s_1)=\{s_1,s_2\}.
\end{aligned}
$$
\end{remark}

Proposition~\ref{descent-inclusion-equality-prop}
immediately implies the following.

\begin{corollary}
\label{even-laced-nonascents-gf}
When $s_0$ is evenly-laced (so $\HatWDes_{R} = \widehat\Des_R$), one has
$$
\begin{aligned}
\sum_{w \in W^+} {\bf t}^{\widehat\Des_{R}(w)} q^{\ell_{R \cup R^{-1}}(w)}
& = \sum_{w \in W^{\{s_0\}}} {\bf t}^{\Des_S(w)} q^{\nu(w)}\\
&= \left[ \sum_{w \in W} {\bf t}^{\Des_S(w)} q_0^{\ell_0(w)} q^{\nu(w)} \right]_{q_0=1,t_0=0}
\end{aligned}
$$
where the elements in $\widehat\Des_{R}(w)$ and $\Des_S(w)$ are identified with their subscripts
as before,
and ${\bf t}^A:=\prod_{j \in A} t_j$.
\end{corollary}

This last generating function for $W$ is easily computed using the techniques
from \cite{Reiner}.

\begin{example} \rm \ \\
\label{second-type-B-example}
Consider the Coxeter system $(W,S)$ of type $B_n$, labelled as in Example~\ref{first-type-B-example}.
Then \cite[\S II, Theorem 3]{Reiner} shows that
$$
\sum_{w \in W} {\bf t}^{\Des_S(w)} q_0^{\ell_0(w)} q^{\nu(w)}
=(-q_0;q)_n [n]!_q \det[ a_{ij} ]_{i,j=-1,0,1,2,\ldots, n-1}
$$
where
$$
a_{ij} =
\begin{cases}
0 & \text{ for }j < i-1 \\
t_i - 1 & \text{ for }j=i-1 \\
\frac{t_i}{(-q_0;q)_{j+1} [j+1]!_q} & \text{ for }  j \geq i= -1 \\
\frac{t_i}{[j-i+1]!_q} & \text{ for } j \geq i \geq 0
\end{cases}
$$
with the convention $t_{-1}=1$.  As an example, for $n=3$, one has
$$
\begin{aligned}
&\sum_{w \in W} {\bf t}^{\Des_S(w)} q_0^{\ell_0(w)} q^{\nu(w)}\\
&\quad =(-q_0;q)_{3} [3]!_q
\det\left[
\begin{matrix}
1 & \frac{1}{(-q_0;q)_1 [1]!)q} &  \frac{1}{(-q_0;q)_2 [2]!)q} &  \frac{1}{(-q_0;q)_3 [3]!)q} \\
t_0-1 & \frac{t_0}{[1]!_q} & \frac{t_0}{[2]!_q} & \frac{t_0}{[3]!_q} \\
0 & t_1-1 & \frac{t_1}{[1]!_q} & \frac{t_1}{[2]!_q}  \\
0 & 0 & t_2-1 & \frac{t_2}{[1]!_q}  \\
\end{matrix}
\right]
\end{aligned}
$$
thus Corollary~\ref{even-laced-nonascents-gf} gives
$$
\begin{aligned}
&\sum_{w \in W^+} {\bf t}^{\HatWDes_{R \cup R^{-1}}(w)} q^{\ell_{R \cup R^{-1}}(w)} \\
&=2(-q;q)_{2} [3]!_q
\det\left[
\begin{matrix}
1 & \frac{1}{2(-q;q)_0 [1]!_q} &  \frac{1}{2(-q;q)_1 [2]!_q} &  \frac{1}{2(-q;q)_2 [3]!_q} \\
-1 & 0 & 0 & 0 \\
0 & t_1-1 & \frac{t_1}{[1]!_q} & \frac{t_1}{[2]!_q}  \\
0 & 0 & t_2-1 & \frac{t_2}{[1]!_q}
\end{matrix}
\right] \\
&= 1  + q (2 t_1 + t_2) + q^2  (3 t_1 + 2 t_2) +
    + q^3  (3 t_1 + t_2 + 2 t_1 t_2)\\
&\qquad \qquad + q^4  (2 t_1 + t_2 + 2 t_1 t_2)
    + q^5  (t_1 + 2 t_1 t_2) + q^6  t_1 t_2.
\end{aligned}
$$
Note that this agrees
with the data in the $1^{st}$ and $6^{th}$ columns from the table
of Section~\ref{large-example} below.
\end{example}

\subsection{Palindromes revisited}
\label{palindromes-revisited}

When $s_0$ is evenly-laced, the set of palindromes for $(W^+,R)$ behaves much more like
the set of reflections in a Coxeter system $(W,S)$, and plays a more closely analogous role.

\begin{theorem}
\label{strong-exchange-property}
Assume $(W,S)$ has $s_0$ evenly-laced.  Then for any $w \in W^+$, one has the following.
\begin{enumerate}
\item[(a)] $\ell_{R \cup R^{-1}} = |\Pal_L(w)|$.
\item[(b)] (Strong exchange property)
For any reduced $(R \cup R^{-1})^*$-word
$$
{\bf r} = (r^{(1)},\ldots,r^{(\nu(w))})
$$
factoring $w$, one has $\Pal_L(w)=\{ p_k \}_{ 1 \leq k \leq \nu(w) }$
where
$$
p_k:=(r^{(1)})^{-1}(r^{(2)})^{-1} \cdots (r^{(k)})^{-1} \cdots (r^{(2)})^{-1} (r^{(1)})^{-1}).
$$
In other words, for a palindrome $p$ and reduced
$(R \cup R^{-1})$-word ${\bf r} = (r^{(1)},\ldots,r^{(\nu(w))})$,
one has
$$
\begin{aligned}
&\ell_{R \cup R^{-1}}(pw) < \ell_{R \cup R^{-1}}(w) \quad \text{ if and only if } \\
&p=p_k \text{ for some }k=1,2,\ldots,\nu(w) \quad \text{ if and only if }\\
&pw = (r^{(1)})^{-1}(r^{(2)})^{-1}\cdots (r^{(k-1)})^{-1} r^{(k+1)} \cdots r^{(\nu(w))}\\
&\qquad \qquad \qquad \text{ for some }k=1,2,\ldots,\nu(w)
\end{aligned}
$$
\item[(c)] The set $\Pal_L(w)$ determines $w$ uniquely.
\end{enumerate}
\end{theorem}
\begin{proof}
Assertions (a) and (b) are immediate from the assertion of equality in
Proposition~\ref{shortening-palindrome-reflection-relation}.

For (c), one must show that for any $w, w' \in W^+$,
if $\Pal_L(w)=\Pal_L(w')$ then $w=w'$.  Via Proposition~\ref{shortening-palindrome-reflection-relation},
it is equivalent to show the following for any $w,w'$ in $W$:
\begin{quote}
$\hat{\T}_L(w)=\hat{\T}_L(w')$ implies $w' \in wW_{\{s_0\}} \, (=\{w,ws_0\})$.
\end{quote}
We will prove this assertion by induction on $\nu(w)$.

In the base case, if $\nu(w)=0$, then $w \in W_{\{s_0\}}$,
which forces $\hat{\T}_L(w')=\hat{\T}_L(w)=\emptyset$, and hence also $w' \in  W_{\{s_0\}}$.

In the inductive step, we make use of the following property \cite[Exercise 1.12]{BB} of $\T_L(w)$:
\begin{equation}
\label{reflection-cocycle}
\begin{aligned}
\T_L(sw) &= \{s\} \,\, \symmdiff \,\, s\T_L(w)s \text{ for any }s \in S,\\
         & \text{ and hence } \\
\hat{\T}_L(s_0w) &=  s_0\hat{\T}_L(w)s_0 \\
\hat{\T}_L(s_iw) &= \{s_i\} \,\, \symmdiff \,\, s_i\hat{\T}_L(w)s_i \text{ for }i=1,2,\ldots,n.
\end{aligned}
\end{equation}
where $A \symmdiff B := (A \setminus B) \sqcup (B \setminus A)$ denotes the symmetric
difference of the sets $A, B$. We treat two cases for $w$.

\noindent
{\bf Case 1}:  $\T_L(w) \cap S \neq \{s_0\}$, say $s_i \in \T_L(w)$ for some  $i=1,2,\ldots,n$.
Then
$$
\hat{\T}_L(s_iw) = \{s_i\} \symmdiff s_i\hat{\T}_L(w)s_i
                  = \{s_i\} \symmdiff s_i\hat{\T}_L(w')s_i
                  = \hat{\T}_L(s_i w').
$$
As $s_i \in \T_L(w)$ implies $\nu(s_iw) < \nu(w)$,
so one can apply induction to conclude that $s_iw' \in s_iw W_{\{s_0\}}$, which implies
$w' \in w W_{\{s_0\}}$ as desired.

\noindent
{\bf Case 2}: $\T_L(w) \cap S = \{s_0\}$.  In this case
$$
\hat{\T}_L(s_0w) = s_0\hat{\T}_L(w)s_0
                  = s_0\hat{\T}_L(w')s_0
                  = \hat{\T}_L(s_0 w'),
$$
and $\nu(s_0w) =\nu(w)$, but $\T_L(s_0w) \cap S \neq \{s_0\}$, so that Case 1 applies.
\end{proof}

Note that we have already seen in Example~\ref{odd-dihedral-example} that,
without the assumption that $s_0$ is evenly-laced, the assertions of
Theorem~\ref{strong-exchange-property} can fail.

%

\subsection{Orders revisited}
\label{orders-revisited}

  When $s_0$ is evenly-laced, the strong exchange property for
palindromes (Theorem~\ref{strong-exchange-property}(b))
has consequences for the weak and strong orders on $W^+$,
analogous to what happens for the weak and strong orders on $W$.

In fact, one can use it to prove the next four propositions,
simply by carrying over the usual proofs from \cite[Chapters 2,3]{BB}, replacing
\begin{enumerate}
\item[$\bullet$]
$S$ with $R \cup R^{-1}$,
\item[$\bullet$]
reflections with palindromes,
\item[$\bullet$]
the usual deletion or strong exchange property with Theorem~\ref{strong-exchange-property}(ii),
\item[$\bullet$]
the trick of writing $w \in W$ as $w=t^2 w$ for a reflection $t \in \T_L(w)$
with the trick of writing $w \in W^+$ as $w=p^{-1} p w$ for a palindrome $p \in  \Pal_L(w)$
(which appeared already in equation~\eqref{rewriting-trick} above).
\end{enumerate}

\begin{proposition}
\label{weak-order-is-palindrome-inclusion}
When $s_0$ is evenly-laced, $u,w \in W^+$ satisfy
$u \leq_{RW} w$ if and only if $\Pal_L(u) \subseteq \Pal_L(w)$.

A similar statement holds for the left weak order $\leq_{LW}$,
replacing left-shortening palindromes $\Pal_L(-)$
with right-shortening palindromes $\Pal_R(-)$.
\end{proposition}

\begin{proposition}
\label{weak-order-is-a-semilattice}
When $s_0$ is evenly-laced, the left, right weak orders on $W^+$ are
meet-semilattices.
\end{proposition}

\begin{proposition}
\label{strong-order-subword-property}
When $s_0$ is evenly-laced, $u,w \in W^+$ satisfy
$u \leq_{LS} w$ if and only if for some (equivalently, every)
reduced $(R \cup R^{-1})^*$-word ${\bf r}=(r^{(1)},\ldots,r^{(\ell)})$
factoring $w$, there exists a reduced $(R \cup R^{-1})^*$-word
factoring $u$ which is a ``subword'' in the following sense:
\begin{quote}
it can be obtained by deleting some of the $r^{(i)}$ from ${\bf r}$ and replacing
any $r^{(i)}$ remaining that have an odd number of letters deleted to their right with
their inverse $(r^{(i)})^{-1}$.
\end{quote}
A similar statement holds for the right strong order $\leq_{RS}$, replacing ``right'' with ``left''.
\end{proposition}

Recall that a poset is {\it thin} if every interval $[x,y]$ of rank $2$ has
exactly $4$ elements, namely $\{x \leq u,v \leq y\}$.

\begin{proposition}
\label{strong-order-is-thin-shellable}
When $s_0$ is evenly-laced, the left, right strong
orders on $W^+$ are thin and shellable, and hence have every open interval
homeomorphic to a sphere.
\end{proposition}

\begin{remark} \rm \ \\
Note that when $s_0$ is not evenly-laced, the strong order need
not be thin, as illustrated by the existence of several upper
intervals of rank $2$ having $5$ elements in
Figure~\ref{Orderings-figure}(b).
\end{remark}

When $(W,S)$ is finite, the examples of $I_2(7), A_3$ from
Figure~\ref{Orderings-figure} show that one need not have a unique
maximum element in any of these orders if $s_0$ is not
evenly-laced. However, if $s_0$ is evenly-laced, there is an
obvious candidate for such a top element, namely $\tau(w_0)$,
where $w_0$ is the longest element of $W$.

\begin{proposition}
\label{finite-case-maximum}
When $(W,S)$ has $s_0$ evenly-laced and $W$ finite, one has
$w_0s_0=s_0w_0$. Furthermore, the element
$\tau(w_0) \in W^+ \cap w_0W_{s_0}$
is the unique maximum element in
all four (left or right, weak or strong) orders on $W^+$.
\end{proposition}
\begin{proof}
For the first assertion note that, by \cite[Proposition 2.3.2]{BB}
$$
\ell_S(w_0 s_0 w_0)
=\ell_S(w_0) - (\ell_S(w_0) - \ell_S(s_0))
= \ell_S(s_0)=1
$$
which shows $w_0 s_0 w_0$ lies in $S$.  But since $w_0^{-1}=w_0$,
it is also conjugate to $s_0$, so in the case
where $s_0$ is evenly-laced, one must have $w_0 s_0 w_0=s_0$, i.e.,
$w_0s_0=s_0w_0$.

To see that $\tau(w_0)$ is the maximum in all four orders,
one can easily check using Proposition~\ref{shortening-palindrome-reflection-relation}
that $\Pal_L(\tau(w_0))=\Pal(W^+)$.
Hence $\tau(w_0)$ is the maximum for the right weak order by
Proposition~\ref{weak-order-is-palindrome-inclusion}.
Since $\tau(w_0)$ is either $w_0$ or $w_0 s_0 = s_0 w_0$,
in either case one has $\tau(w_0)^{-1}=\tau(w_0)$, and hence
it is also the maximum for the left weak order.
It is then also the maximum for the left and right strong orders because they are stronger than
the corresponding weak orders.
\end{proof}

\section{The case of a leaf node}
\label{leaf-section}

  The presentation \eqref{general-W-plus-presentation} for $W^+$
becomes very close to a Coxeter presentation when $s_0$ is a {\it leaf} node, that is,
$s_0$ commutes with $s_2,\ldots,s_n$, i.e., one has $m_{0i}=2$ for $i=2,\ldots,n$
(although $m_{01}$ may be greater than $2$).
Note that every (irreducible) finite and affine Coxeter system $(W,S)$, with the exception of
the family  $\tilde{A}_n$, has Coxeter diagram shaped like a tree, and hence will have some leaf
node $s_0$.

\subsection{Nearly Coxeter presentations}

\begin{proposition}
\label{leaf-presentation-prop}
Let $(W,S)$ be a Coxeter system with $S=\{s_0,s_1,\ldots,s_n\}$ and
$s_0$ a leaf node.
Then $W^+$ is generated by the set
$$
R:=\{r_i=s_0s_i\ |\ s_i\in S\setminus s_0 \}
$$
with the following presentation:
\begin{equation}
\label{leaf-presentation}
\begin{aligned}
W^+ \cong \langle R=\{r_1,\ldots,r_n\} :
r_1^{m_{01}}=r_i^2 = (r_ir_j)^{m_{ij}}=e \text{ for }1 \leq i < j \leq n \rangle,
\end{aligned}
\end{equation}
where $m_{ij}$ is the order of $s_is_j$ and $s_1$ is the neighbor
of the leaf $s_0$.
\end{proposition}
\begin{proof}
Starting with the presentation in  \eqref{general-W-plus-presentation},
note that given any $1 \leq i < j \leq n$, the relation
$(r_i^{-1}r_j)^{m_{ij}}=e$ is equivalent to
$(r_j^{-1}r_i)^{m_{ij}}=e$ by taking the inverse of both sides.
However, since $j \geq 2$, one has $r_j^2=e$ and so $r_j^{-1}=r_j$.
Thus this relation is equivalent to $(r_jr_i)^{m_{ij}}=e$, which is
also equivalent to $(r_i r_j)^{m_{ij}}=e$ via conjugation by $r_j$.
\end{proof}

\begin{defn} \rm \ \\
Call a presentation for an abstract group having
the form in \eqref{leaf-presentation} a {\it nearly Coxeter}
presentation, meaning that all but one of the generators $r_i$ is
an involution and all other relations are of the form
$(r_i r_j)^{m_{ij}}$ for some $m_{ij} \in \{2,3,4,\ldots\} \cup
\{\infty\}$.
\end{defn}

\begin{corollary}
Every abstract group $A$ with a nearly Coxeter presentation is isomorphic to the
alternating subgroup $W^+$ of some Coxeter system $(W,S)$.

In particular, if
$A$ is finite and has a nearly Coxeter presentation, then it is isomorphic to
a product
\begin{equation}
\label{A-decomposition}
A \cong W_0^+ \times W_1 \times \cdots \times W_r
\end{equation}
in which each of the $(W_i,S_i)$ are finite irreducible Coxeter systems
(and hence classified).
\end{corollary}
\begin{proof}
If $A$ is an abstract group with a nearly Coxeter presentation, as in \eqref{leaf-presentation},
one can write down a corresponding Coxeter system $(W,S)$ as in \eqref{W-presentation}.
Theorem~\ref{general-W-plus-presentation} then shows that $A \cong W^+$.

Furthermore, if $A$ is finite, then since $A \cong W^+$ and $[W:W^+]=2$, one concludes
that $W$ is also finite.  Consequently
$$
W \cong  W_0 \times W_1 \times \cdots \times W_r
$$
for some finite irreducible Coxeter systems $(W_i,S_i)$.  Without loss of generality,
one can index so that $s_0,s_1$ belong to $(W_0,S_0)$.  The
isomorphism \eqref{A-decomposition} then follows from examining the presentation.
\end{proof}

\section{The case of an even leaf node}
\label{even-leaf-section}

  When the distinguished node $s_0$ is both a leaf and evenly-laced,
that is, $m_{01}$ is even and $m_{0j}=2$ for $j=2,3,\ldots,n$, we
shall say that $s_0$ is an {\it even leaf}.  In this situation
$(W^+,R)$ has an amazingly close connection to
the index $2$ subgroup $W':=\ker \chi_0$ of $W$, which will turn
out to have a Coxeter structure $(W',S')$ of its own. Note that in
every finite and affine Coxeter system containing an evenly-laced
node $s_0$, namely types $B_n(=C_n), \tilde{B}_n, \tilde{C}_n$,
this evenly-laced node is actually an even leaf\footnote{with a
single affine exception: the affine type $\tilde{C}_2$ has the
middle node in its diagram evenly-laced, but not an even leaf!},
to which the results below apply\footnote{ While the combinatorics
of $W^+$ and $W'$ seems to be similar, the combinatorics of other
subgroups of index 2 seems to be different; in particular, no
nearly Coxeter presentation for these groups is known; see, e.g.,
\cite{Bernstein}.}.

\subsection{The Coxeter system $(W',S')$} \ \\
\label{W'-presentation-section}

Assume $(W,S)$ is a Coxeter system with $S=\{s_0,s_1,\ldots,s_n\}$ having $s_0$ as an
even leaf.  Since $s_0$ is
evenly laced, recall that one has the linear character
$\chi_0: W \rightarrow \{\pm 1\}$ from \eqref{chi-zero-definition},
taking value $-1$ on $s_0$ and $+1$ on all other $s_j \in S$.
Let $W':=\ker \chi_0$, a subgroup of $W$ of index $2$.

We wish to show that $W'$ is a reflection subgroup of $W$, and has a natural
Coxeter presentation $(W',S')$ extremely close to $(W,S)$.
Let $S':=\{t_1,t_2,\ldots,t_n\} \cup \{t_1'\}$ be a set, and
consider the set map
$$
\begin{matrix}
S' & \overset{f}{\longrightarrow} &W' &\\
t_j& \overset{f}{\longmapsto} & s_j   &\text{ for }j=1,2\ldots,n \\
t_1'&\overset{f}{\longmapsto} & s_0 s_1 s_0&
\end{matrix}.
$$

\begin{proposition}
\label{W'-Coxeter-presentation}
The set map $f$ above extends to an isomorphism
\begin{equation}\label{eq-W'-Coxeter-presentation}
\begin{aligned}
W' \cong
&\langle S'=\{t_1,\ldots,t_n\} \cup \{t_1'\} : \\
&\qquad \qquad (t_i)^2=(t_1')^2=(t_it_j)^{m_{ij}}=e \text{ for }1\leq i\leq j \leq n,\\
&\qquad \qquad (t_1' t_j)^{m_{1j}}=e,\\
&\qquad \qquad (t_1' t_1)^{\frac{m_{01}}{2}}=e \rangle.
\end{aligned}
\end{equation}
which makes $(W',S')$ a Coxeter system.
\end{proposition}

A schematic picture of the relation between the Coxeter diagrams of $(W,S)$
and $(W',S')$ was shown in Figure~\ref{Oriflamme-figure}.  Note that the embedding $W' \subset W$
as a reflection (but not parabolic) subgroup generalizes the finite/affine Weyl
group inclusions
$$
\begin{aligned}
W(D_n) &\subseteq W(B_n)(=W(C_n)), \\
W(\tilde{D}_n) &\subseteq W(\tilde{B}_n), \text{ and }\\
W(\tilde{B}_n) & \subseteq W(\tilde{C}_n)
\end{aligned}
$$
in which one always has $m_{01}=4$ so that $t_1, t_1'$ commute, and are a
pair of {\it oriflamme/fork} nodes at the end of the Coxeter diagram for $(W',S')$.

\begin{proof} (of Proposition~\ref{W'-Coxeter-presentation}).

We employ a similar trick to Bourbaki's from
Proposition~\ref{Bourbaki-exercise}. Consider the abstract group
$G$ with the Coxter presentation given on the right side of
\eqref{eq-W'-Coxeter-presentation}. Since $t_1,t_1'$ play
identical roles in this presentation, the set map $\beta: S'
\rightarrow G$ which fixes $t_2,\ldots,t_n$ and swaps $t_1', t_1$
extends to an involutive group automorphism $\beta: G \rightarrow
G$.

Thus the group $\ZZ/2\ZZ=\{1,\beta\}$ acts on $G$, and one can form
the semidirect product $G \rtimes \ZZ/2\ZZ$ in which
$(g_1\beta^i) \cdot (g_2 \beta^j )=g_1 \beta^i(g_2) \cdot \beta^{i+j}$.
This has the following presentation:
$$
\begin{aligned}
G \rtimes \ZZ/2\ZZ \cong
&\langle t_1,\ldots,t_n, \beta: \\
&\qquad \qquad \beta^2=(t_1')^2=(t_it_j)^{m_{ij}}=e \text{ for }1\leq i\leq j \leq n,\\
&\qquad \qquad (t_1' t_j)^{m_{1j}}=e,\\
&\qquad \qquad (t_1' t_1)^{\frac{m_{01}}{2}}=e, \\
&\qquad \qquad \beta t_j = t_j \beta \text{ for }2 \leq j \leq n,\\
&\qquad \qquad \beta t_1 = t_1' \beta \rangle.
\end{aligned}
$$

We claim that the following maps $g, f$
are well-defined and inverse isomorphisms:
$$
\begin{matrix}
W       & \overset{g}{\longrightarrow}  & G \rtimes \ZZ/2\ZZ &\\
   s_i  & \longmapsto     & t_i   & \text{ for }i=1,\ldots,n\\
   s_0  & \longmapsto     & \beta &\\
     &  &  & \\
G \rtimes \ZZ/2\ZZ & \overset{f}{\longrightarrow} & W & \\
            t_i  & \longmapsto     & s_i &\text{ for }i=1,\ldots,n \\
            t_1'  & \longmapsto     & s_0 s_1 s_0 & \\
          \beta  & \longmapsto     & s_0 &
\end{matrix}
$$
Here are the relations in $(W,S)$ going to relations in $G \rtimes \ZZ/2 \ZZ$ needed
to check that $f$ is well-defined:
$$
\begin{aligned}
s_0^2=e & \quad \mapsto \quad  \beta^2=e \\
(s_i s_j)^{m_{ij}}=e & \quad  \mapsto  \quad (t_i t_j)^{m_{ij}}=e \text{ for }1 \leq i \leq j \leq n\\
(s_0 s_j)^2=e &  \quad \mapsto  \quad (\beta t_j)^2=\beta t_j \beta t_j = \beta \beta t_j t_j = e \text{ for }2 \leq j \leq n \\
(s_0 s_1)^{m_{01}}=e &  \quad \mapsto  \quad (\beta t_1)^{m_{01}}
                      = \underbrace{ \beta t_1 \cdot \beta t_1 \cdots }_{m_{01}\text{ times}}
                      = \underbrace{ \beta t_1 t_1' \beta \cdot
                           \beta t_1 t_1' \beta  \cdots }_{\frac{m_{01}}{2}\text{ times}} \\
                     & \qquad \qquad \qquad \qquad = \beta ( t_1 t_1' )^{\frac{m_{01}}{2}} \beta =e.
\end{aligned}
$$
Here are the relations in $G \rtimes \ZZ/2 \ZZ$ going to relations in $(W,S)$ needed to
check that $g$ is well-defined:
$$
\begin{aligned}
\beta^2=e &  \quad \mapsto  \quad s_0^2=e \\
(t_i t_j)^{m_{ij}}=e &  \quad \mapsto  \quad (s_i s_j)^{m_{ij}}=e \text{ for }1 \leq i \leq j \leq n\\
(t_1')^2=e &  \quad \mapsto  \quad (s_0 s_1 s_0)^2 = e \\
(t_1' t_j)^{m_{1j}}=e &  \quad \mapsto  \quad (s_0 s_1 s_0 s_j)^{m_{1j}} = (s_0 s_1 s_j s_0)^{m_{1j}} = s_0 (s_1 s_j)^{m_{1j}} s_0 = e\\
(t_1 t_1')^{\frac{m_{01}}{2}}=e &  \quad \mapsto  \quad (s_1 s_0 s_1
s_0)^{\frac{m_{01}}{2}} = (s_1 s_0)^{m_{01}} = e\\
\beta t_j = t_j \beta &  \quad \mapsto  \quad s_0 s_j=s_j s_0 \text{ for } 2\le
j\le n\\
\beta t_1 = t_1' \beta &  \quad \mapsto  \quad s_0 s_1 = s_0 s_1 s_0 s_0.
\end{aligned}
$$

\noindent
Once one knows that $f, g$ are well-defined, it is easily checked that they
are inverse isomorphisms by checking this on generators.

Since $f(G) \subseteq W'$, and both $W', f(G)$
are subgroups of $W$ of index $2$, it must be that $f(G)=W'$.  Hence $f$
restricts to the desired isomorphism presenting $W'$ as the Coxeter group $G$.
\end{proof}

\subsection{Relating $W^+$ to $W'$}
\label{W-plus-to-W'-prime-section}

  We next discuss the tight relation between $(W^+,R)$ and $(W',S')$, which
is mediated by the following map.

\begin{proposition}
When $s_0$ is an even leaf in $(W,S)$, the following formulae
$$
\begin{matrix}
W^+ &\overset{\theta}{\longrightarrow} & W' &\\
  w & \longmapsto& w \cdot s_0^{\ell_{R \cup R^{-1}}(w)} = &
   \begin{cases}
      w    & \text{ if } w \in W' \\
      ws_0 & \text{ if } w \not\in W'
   \end{cases}
\end{matrix}
$$
define the same set map $\theta:W^+ \rightarrow W'$.
In other words, $\theta(w)$ is the unique element in the coset
$wW_{\{s_0\}}=\{w,ws_0\}$ that lies in $W'$.

Furthermore, $\theta$ is a bijection, and equivariant for the action of the subgroup
$W^+ \cap W'$ by left-multiplication on $W^+$ and $W'$.
\end{proposition}

\begin{proof}
Note that
$$
\chi_0(r_i)=\chi_0(s_0 s_i) = -1 = \chi_0(s_i s_0) = \chi_0(r^{-1}_i) \text{ for all }i
$$
and hence $\chi_0(w)=(-1)^{\ell_{R \cup R^{-1}}(w)}$.  This shows the equivalence of the two formulae
for $\theta(w)$.

The $W^+ \cap W'$-equivariance of $\theta$ follows from either formula.
Bijectivity of $\theta$ follows, for example,
since one can check that the map $\tau: W \rightarrow W^+$
from Definition~\ref{tau-definition} when restricted to $W'$
satisfies $\tau|_{W'}=\theta^{-1}$.
\end{proof}

Note that the bijection $\theta: W^+ \rightarrow W'$ is {\it not} a group isomorphism,
and that $W^+, W'$ are generally not isomorphic as groups.  For example, when $(W,S)$ is a dihedral
Coxeter system $I_2(m)$ with $m$ even, $W^+$ is always cyclic of order $m$, while $W'$ is not
cyclic for $m \geq 4$.

Nevertheless, the map $\theta$ is about as close as one can get to an
isomorphism of the presentations $(W^+,R)$ and $(W',S')$, in that $\theta$
lifts to the following map on words in the generating sets.

\begin{defn} \rm \ \\
When $s_0$ is an even leaf in $(W,S)$, define a set map $\Theta: (R \cup R^{-1})^* \rightarrow (S')^*$
by mapping a word ${\bf r}=(r^{(1)},\ldots,r^{(\ell)})$ one letter at a time according to the
following rules:
$$
\begin{matrix}
  r_j & \longmapsto & t_j  \quad \text{ for }j=2,3,\ldots,n \\
   &  &    \\
  r_1 & \longmapsto&
  \begin{cases}
      t_1    & \text{ if } r_1 = r^{(k)} \text{ with }k\text{ even }, \\
      t'_1    & \text{ if } r_1 = r^{(k)} \text{ with }k\text{ odd }, \\
   \end{cases}  \\
   &  &   \\
  r^{-1}_1 & \longmapsto&
  \begin{cases}
      t'_1    & \text{ if } r_1 = r^{(k)} \text{ with }k\text{ even }, \\
      t_1    & \text{ if } r_1 = r^{(k)} \text{ with }k\text{ odd }. \\
   \end{cases}  \\
\end{matrix}
$$
\end{defn}

  The maps $\theta, \Theta$ are related as follows.

\begin{proposition}
\label{Theta-lifts-theta}
Let  $s_0$ be an even leaf in $(W,S)$.  Then for any $(R \cup R^{-1})^*$-word ${\bf r}$
of length $\ell$ that factors $w \in W^+$, its image $\Theta({\bf r})$ is an
$S^*$-word of the same length that factors $\theta(w)$.
\end{proposition}
\begin{proof}
Given ${\bf r}=(r^{(1)},\ldots,r^{(\ell-1)},r^{(\ell)})$ of length
$\ell$ factoring $w \in W^+$, denote by $w'$ the element in $W'$
factored by $\Theta({\bf r})$.  One must show that $w^{-1} w' =
s_0^\ell$, where $\ell = \ell_{R \cup R^{-1}}(w)$.

Proceed by induction on $\ell$, where the base case $\ell=0$ is trivial.
In the inductive step, let $u$ denote the element in $W^+$ factored by $(r^{(1)},\ldots,r^{(\ell-1)})$,
and $u'$ the element in $W'$ factored by $\Theta((r^{(1)},\ldots,r^{(\ell-1)})$.
By induction, $u^{-1} u' = s_0^{\ell-1}$. Then
$$
\begin{aligned}
w^{-1} w' &= (r^{(\ell)})^{-1} \cdot u^{-1} u' \cdot \Theta( r^{(\ell)} )  \\
         & = (r^{(\ell)} )^{-1} \cdot s_0^{\ell-1} \cdot \Theta( r^{(\ell)} ).
\end{aligned}
$$
which one must show coincides with $s_0^\ell$.  Consider the following cases:

\noindent {\bf Case 1.} $r^{(\ell)} = r_j$ for some
$j=2,3,\ldots,n$. Then $r^{(\ell)} = s_0 s_j$ and $\Theta(
r^{(\ell)} ) = t_j = s_j$, so one gets
$$
(s_0 s_j)^{-1} \cdot s_0^{\ell-1} \cdot s_j = s_0^\ell
$$
because $s_0, s_j$ commute.

\noindent
{\bf Case 2a.} $r^{(\ell)} = r_1$  and $\ell$ is even.
Then $r^{(\ell)} = s_0 s_1$ and $\Theta( r^{(\ell)} ) = t_1 = s_1$, so one gets
$$
(s_0 s_1)^{-1} \cdot s_0^{\ell-1} \cdot s_1
= s_1 s_0 \cdot s_0^{\ell-1} \cdot s_1
= s_1 s_0^\ell s_1 = e = s_0^{\ell}.
$$

\noindent
{\bf Case 2b.} $r^{(\ell)} = r_1$  and $\ell$ is odd.
Then $r^{(\ell)} = s_0 s_1$ and $\Theta( r^{(\ell)} ) = t'_1 = s_0 s_1 s_0$, so one gets
$$
(s_0 s_1)^{-1} \cdot s_0^{\ell-1} \cdot s_0 s_1 s_0
= s_1 s_0 \cdot s_0^{\ell-1} \cdot s_0 s_1 s_0
= s_0 = s_0^{\ell}.
$$

\noindent {\bf Case 3a.} $r^{(\ell)} = r^{-1}_1$  and $\ell$ is
even. Then $r^{(\ell)} = s_1 s_0$ and $\Theta( r^{(\ell)} ) = t'_1
= s_0 s_1 s_0$, so one gets
$$
(s_1 s_0)^{-1} \cdot s_0^{\ell-1} \cdot s_0 s_1 s_0
= s_0 s_1 \cdot s_0^{\ell-1} \cdot s_0 s_1 s_0
= e = s_0^{\ell}.
$$

\noindent
{\bf Case 3b.} $r^{(\ell)} = r^{-1}_1$  and $\ell$ is odd.
Then $r^{(\ell)} = s_1 s_0$ and $\Theta( r^{(\ell)} ) = t_1 = s_1$, so one gets
$$
(s_1 s_0)^{-1} \cdot s_0^{\ell-1} \cdot s_1
= s_0 s_1 \cdot s_0^{\ell-1} \cdot s_1
= s_0 = s_0^{\ell}.
$$
\end{proof}

\begin{corollary}
\label{tight-relation}
Let $s_0$ be an even leaf in $(W,S)$.  Then for any $w \in W^+$, the bijections $\theta, \Theta$
have the following properties:
\begin{enumerate}
\item[(i)] $\ell_{R \cup R^{-1}}(w) = \ell_{S'}(\theta(w))$.
\item[(ii)] $\Theta$ bijects the set of reduced $(R \cup R^{-1})^*$-words for $w$
with the reduced $(S')^*$-words for $\theta(w)$
\item[(iii)]  Given any $w \in W^+$, the bijection $R \cup R^{-1} \rightarrow S'$ defined by
$$
\begin{aligned}
r_j          &\mapsto t_j \text{ for }j=2,3,\ldots,n \\
r_1,r_1^{-1} &\mapsto
\begin{cases}
t'_1,t_1 &\text{ if }\ell_{R \cup R^{-1}}(w) \text{ is even},\\
t_1,t'_1 &\text{ if }\ell_{R \cup R^{-1}}(w) \text{ is odd},\\
\end{cases}
\end{aligned}
$$
bijects $\Des_{R \cup R^{-1}}(w)$ with $\Des_{S'}(\theta(w))$.
\item[(iv)] The map $\theta$ is a poset isomorphism $(W^+, \leq_{RW}) \rightarrow (W', \leq_{RW} )$.
\item[(v)] The map $\theta$ is a poset isomorphism $(W^+, \leq_{RS}) \rightarrow (W', \leq_{S} )$.
\end{enumerate}
\end{corollary}
\begin{proof}
Assertions (i),(ii),(iii) and (iv) are straightforward from
Proposition~\ref{Theta-lifts-theta}, while assertion (v) follows
from it via 
Proposition~\ref{strong-order-subword-property}.
\end{proof}

\begin{corollary}
\label{tight-gf-relation}
When $s_0$ is an even leaf in $(W,S)$, one has
\begin{equation}
\label{even-leaf-refined-gf}
\sum_{w \in W^+} {\bf t}^{\Des_{R \cup R^{-1}}(w)} q^{\ell_{R \cup R^{-1}}(w)}
  = {\bf \Theta}\left[ \sum_{w \in W'} {\bf t}^{\Des_{S'}(w)} q^{\ell_{S'}(w)} \right].
\end{equation}
where ${\bf \Theta}$ is the obvious operator on monomials ${\bf
t}^A q^\ell$ corresponding to the mapping from
Corollary~\ref{tight-relation}(iii).  In particular, letting
$\des_{R \cup R^{-1}}(w):=|\Des_{R \cup R^{-1}}(w)|$,
\begin{equation}
\label{even-leaf-des-length-gf}
\sum_{w \in W^+} t^{\des_{R \cup R^{-1}}(w)} q^{\ell_{R \cup R^{-1}}(w)}
  =  \sum_{w \in W'} t^{\des_{S'}(w)} q^{\ell_{S'}(w)},
\end{equation}
and
\begin{equation}
\label{even-leaf-length-gf}
W^+(R \cup R^{-1};q) =W'(S';q).
\end{equation}
\end{corollary}

\begin{example} \rm \ \\
\label{type-B-length-gf-again}
Let $(W,S)$ be the Coxeter system of type $B_n$, labelled as in Example~\ref{first-type-B-example}.
Here $(W',S')$ is the Coxeter system of type $D_n$, whose exponents are known to be
$n-1,1,3,5,\ldots,2n-5,2n-3$.  Hence one can rederive the length generating function
for $(W^+,R)$ using Corollary~\ref{tight-gf-relation} and a well known result in
the theory of Coxeter groups
(see, e.g., \cite[Theorem 7.1.5]{BB} or \cite[Theorem 3.15]{Humphreys}) as follows:
$$
W^+(R \cup R^{-1};q)
 = W'(S';q)
 = [n]_q \prod_{j=1}^{n-1} [2j]_q.
$$
\noindent Furthermore, \cite[Theorem 7]{Reiner} gives generating
functions incorporating the distributions of descents and length
simultaneously for all groups $D_n$, and hence
equation~\eqref{even-leaf-des-length-gf} allows one to derive the
generating functions of $\des_{R \cup R^{-1}}$ and $\ell_{R \cup
R^{-1}}$ simultaneously for $W^+$ of all of the groups $B_n$ when
$s_0$ is chosen to be the even leaf.
When $n=3$ this gives, for example
$$
\begin{aligned}
\sum_{w \in W^+} t^{\des_{R \cup R^{-1}}(w)} q^{\ell_{R \cup R^{-1}}(w)}
 &= 1 + 3 q t + q^2  (4 t + t^2 ) + q^3  (3 t + 3 t^2 )  \\
 &\qquad \qquad \qquad +  q^4  (t + 4 t^2 ) + 3 q^5  t^2 + q^6  t^3,
\end{aligned}
$$
which agrees with the data in the $1^{st}$ and $5^{th}$ columns of the
table in Section~\ref{large-example} below.
\end{example}

\begin{example} \rm \ \\
\label{affine-length-gfs} For $(W,S)$ of affine type
$\tilde{C_n}$, one has $(W',S')$ equal to the affine Coxeter
system of type $\tilde{B}_n$. Using
Corollary~\ref{tight-gf-relation}, the known exponents
$1,3,5,\ldots,2n-1$ for the finite type $B_n$, and Bott's formula
(see, e.g., \cite[Theorem 7.1.10]{BB} or \cite[\S 8.9]{Humphreys})
for the length generating function of an
affine Weyl group, one has that
$$
W^+(R \cup R^{-1};q)
 = W'(S';q)
 = \prod_{j=1}^{n}\frac{[2j]_q}{1-q^{2j-1}}.
$$

Similarly, for $(W,S)$ of affine type $\tilde{B_n}$,
one has that $(W',S')$ is the affine Coxeter system of type $\tilde{D}_n$,
and one derives
$$
W^+(R \cup R^{-1};q)
 =  W'(S';q)
 = \frac{ [n]_q }{1-q^{n-1}} \prod_{j=1}^{n-1}\frac{[2j]_q}{1-q^{2j-1}}.
$$
\noindent A refinement may be obtained using~\cite[Theorems 7 and
8]{Reiner}, which give generating functions incorporating the
distributions of descents and length simultaneously for all groups
$\tilde{B}_n, \tilde{D}_n$. Hence,
equation~\eqref{even-leaf-des-length-gf} allows one to derive the
generating functions of $\des_{R \cup R^{-1}}$ and $\ell_{R \cup
R^{-1}}$ simultaneously for $W^+$ of all of the
groups $\tilde{C}_n, \tilde{B}_n$, when $s_0$ is chosen to be an even leaf.
\end{example}


\subsection{The example of $B_3$}
\label{large-example}
  We compute here $W^+$ in a reasonably large example with $s_0$ an even leaf,
as an illustration of some of the preceding results.

  Consider again the Coxeter system $(W,S)$ of type $B_n$, labelled
as in Example~\ref{first-type-B-example}. The Coxeter system
$(W',S')$ is of type $D_3 (=A_3)$, with $S'=\{t_1',t_1,t_2\}$
having $m(t_1,t_1')=2, m(t_1,t_2)=m(t_1',t_2)=3$. In the table
below, the first three columns give the elements $w$ of $W'$
according to their $S'$-length, giving the list of $S'$-reduced
words for each (abbreviating $t_1,t_1',t_2$ by $1,1',2$) and their
descent set $\Des_{S'}(w)$.  The remaining  columns give the
corresponding element of $W^+$ with its $(R \cup
R^{-1})^*$-reduced words (abbreviating $r_1,r_1^{-1},r_2$ by
$1,\bar{1},2$), its nonascent set $\WDes_{R\cup R^{-1}}(=\Des_{R \cup R^{-1}})$
and its symmetrized nonascent set $\HatWDes_R(=\widehat\Des_R)$.

%
$$
\begin{matrix}
 \ell_{S'} &
   (S')^* &
     \Des_{S'} &
             (R \cup R^{-1})^*&
                 \WDes_{R\cup R^{-1}} & \HatWDes_R \\
=\ell_{R \cup R^{-1}} &
       \text{reduced words} &
             &
       \text{reduced words} & =\Des_{R \cup R^{-1}}& =\widehat\Des_R \\
    &  &   &  &   \\
 0  & \emptyset    & \emptyset  & \emptyset & \emptyset & \emptyset\\
    &  &   &  &   \\
    &  &   &  &   \\
 1  & 1  & 1   & \bar{1} &  1  & 1\\
    &  &   &  &   \\
    & 1' & 1'   & 1      & \bar{1} & 1 \\
    &  &   &  &   \\
    & 2  & 2   & 2      & 2 & 2\\
    &  &   &  &   \\
    &  &   &  &   \\
 2  & 12  & 2   & \bar{1}2 & 2 & 2  \\
    &  &   &  &   \\
    & 1'2 & 2   & 12      & 2 & 2\\
    &  &   &  &   \\
    & 21  & 1   & 21      & \bar{1} & 1 \\
    &  &   &  &   \\
    & 21'  & 1'   & 2\bar{1}      & 1 & 1 \\
    &  &   &  &   \\
    & 11',1'1  & 1,1'   & \bar{1}\bar{1},11      & 1,\bar{1} & 1 \\
    &  &   &  &   \\
    &  &   &  &   \\
 3  & 121,212  & 1,2   & \bar{1}2\bar{1},212 &  1,2  & 1,2 \\
    &  &   &  &   \\
    & 1'21',21'2 & 1',2   & 121,2\bar{1}2 & \bar{1},2 & 1,2 \\
    &  &   &  &   \\
    & 11'2,1'12  & 2   & \bar{1}\bar{1}2,112 & 2 & 2 \\
    &  &   &  &   \\
    & 211',21'1  & 1,1'   & 211,2\bar{1}\bar{1} & 1,\bar{1} & 1 \\
    &  &   &  &   \\
\end{matrix}
$$

$$
\begin{matrix}
 \ell_{S'} &
   (S')^* &
     \Des_{S'} &
             (R \cup R^{-1})^*&
                 \WDes_{R\cup R^{-1}}  & \HatWDes_R \\
=\ell_{R \cup R^{-1}} &
       \text{reduced words} &
             &
       \text{reduced words} & =\Des_{R \cup R^{-1}}& =\widehat\Des_R\\
    &  &   &  &   \\
    & 121'  & 1'   & \bar{1}21      & \bar{1} & 1 \\
    &  &   &  &   \\
    & 1'21  & 1   & 12\bar{1} & 1 & 1 \\
    &  &   &  &   \\
    &  &   &  &   \\
 4  & 121'1, 1211',2121' & 1,1'   & \bar{1}211, \bar{1}2\bar{1}\bar{1},212\bar{1}
            &  1,\bar{1}  & 1 \\
    &  &   &  &   \\
    & 1'211',1'21'1,& 1,1'   & 12\bar{1}\bar{1},1211,
            & 1,\bar{1} & 1 \\
    & 21'21  &   &2\bar{1}21  &   \\
    &  &   &  &   \\
    & 11'21,1'121, & 1,2   & \bar{1}\bar{1}21,1121,
            & \bar{1}, 2 & 1,2 \\
    &1'212  &   &12\bar{1}2  &   \\
    &  &   &  &   \\
    & 1'121', 11'21' & 1',2
         & 112\bar{1},\bar{1}\bar{1}2\bar{1},     & 1,2 & 1,2 \\
    & 121'2 &   & \bar{1}212  &   \\
    &  &   &  &   \\
    & 211'2,21'12  & 2   & 2112,2\bar{1}\bar{1}2  & 2 & 2 \\
    &  &   &  &   \\
    &  &   &  &   \\

%
5  & 121'12, 1211'2,   & 1',2      & \bar{1}2112,
\bar{1}2\bar{1}\bar{1}2,
    &  \bar{1},2 & 1,2  \\
    & 2121'2, 211'21',  &   & 212\bar{1}2, 21121, &  \\
    & 21'121' &             &  2\bar{1}\bar{1}21 &   \\
    &  &   &  &   \\
    & 211'21, 21'121,   & 1,2
    & 2112\bar{1}, 2\bar{1}\bar{1}2\bar{1}, &  1,2 & 1,2   \\
    & 21'212, 1'21'12,  &   & 2\bar{1}212, 12112,  &  \\
    &  1'211'2          &   & 12\bar{1}\bar{1}2 &   \\
    &  &   &  &   \\
    & 121'21, 11'21'1, & 1,1'  & \bar{1}212\bar{1}, \bar{1}\bar{1}2\bar{1}\bar{1},  &  1,\bar{1} & 1  \\
    & 1'121'1, 11'211', &   &  112\bar{1}\bar{1}, \bar{1}\bar{1}211, & \\
    & 1'1211' , 1'2121' &   &   11211,12\bar{1}21 &   \\
      &  &   &  &   \\
      &  &   &  &   \\
 6  & 1211'21,121'121, & 1,1',2 & \bar{1}2\bar{1}\bar{1}21,\bar{1}21121, &  1, \bar{1},2 & 1,2  \\
    & 2121'21, 121'212, & & 212\bar{1}21, \bar{1}212\bar{1}2, & \\
    & 211'21'1, 11'21'12, & & 211211, \bar{1}\bar{1}2\bar{1}\bar{1}2, & \\
    & 21'121'1, 11'211'2, & & 2\bar{1}\bar{1}211, \bar{1}\bar{1}2112, & \\
    & 211'211', 1'121'12, & & 2112\bar{1}\bar{1}, 112\bar{1}\bar{1}2, & \\
    & 21'11'11', 1'1211'2, & & 2\bar{1}\bar{1}2\bar{1}\bar{1}, 112112, & \\
    & 21'2121', 1'2121'2, & & 2\bar{1}212\bar{1}, 12\bar{1}212, & \\
    & 1'21'121', 1'211'21', & & 12112\bar{1}, 12\bar{1}\bar{1}2\bar{1}, & \\
\end{matrix}
$$


\section{Acknowledgements}
The authors thank John Stembridge for pointing them to the Bourbaki exercise
from Section~\ref{Bourbaki-subsection} at an early stage of this investigation.


\end{document}